# Tomographic projection optimization for volumetric additive manufacturing with general band constraint Lp-norm minimization

Chi Chung Li[1*], Joseph Toombs[1], Hayden K. Taylor[1], Thomas J. Wallin[2†]

[1]Department of Mechanical Engineering, University of California, Berkeley, Etcheverry Hall, 2521 Hearst Ave, Berkeley, CA 94720, USA

[2]Department of Materials Science and Engineering, Massachusetts Institute of Technology, 6-113, 77 Massachusetts Ave., Cambridge, MA 02139, USA

[*]Corresponding author: ccli@berkeley.edu

[†]Corresponding author: tjwallin@mit.edu

## Abstract

Tomographic volumetric additive manufacturing is a rapidly growing fabrication technology that enables rapid production of 3D objects through a single build step. In this process, the design of projections directly impacts geometric resolution, material properties, and manufacturing yield of the final printed part. Herein, we identify the hidden equivalent operations of three major existing projection optimization schemes and reformulate them into a general loss function where the optimization behavior can be systematically studied, and unique capabilities of the individual schemes can coalesce. The loss function formulation proposed in this study unified the optimization for binary and greyscale targets and generalized problem relaxation strategies with local tolerancing and weighting. Additionally, this formulation offers control on error sparsity and consistent dose response mapping throughout initialization, optimization, and evaluation. A parameter-sweep analysis in this study guides users in tuning optimization parameters for application-specific goals.

## 1 Introduction

Tomographic volumetric additive manufacturing (VAM) [1–3] construct 3D objects by irradiating a rotating vat of photopolymer resin with a 2D dynamic light pattern. The cumulative photodosage in a given volume dictates the extent of photochemical reactions and thereby local material response. Above a threshold reaction conversion, the material polymerizes into a free-standing solid. VAM offers numerous benefits such as high fabrication speed (up to 4500 $mm^3 \cdot min^{-1}$ [4]), smooth part surface (with roughness $R_a$ down to 6 $nm$ [5]), support-free printing of overhang features, the ability to overprint onto existing structure and internal patterning in preassembled material. These fabrication advantages are well suited for application involving printing of optical components [5,6], shape memory polymers [7], object-embedded organogel [8], biocompatible hydrogels [9,10], multi-material constructs [11] and for printing in microgravity environments [4]. In practice, the geometry, surface finish and internal conversion profile greatly impact the performance and function of the printed part. In VAM, these factors are influenced by both the tomographic exposure step and post-processing steps. In particular, the exposure step controls local material response and dictates the spatial conversion profile of the material in the as-patterned state.

In the subsequent postprocessing steps, the as-patterned material goes through development steps to remove the low-conversion and ungelled phase, and optionally further post-cured by optical flood exposure and thermal treatment. Since the postprocessing steps only provide global transformations, spatial control of conversion and other associated properties crucially depends on the tomographic projections used in the initial exposure step. The goal of this work is to identify and define a general mathematical structure of the projection optimization problem such that these local material responses can be more finely controlled.

The tomographic projection optimization problem in VAM is closely related to the inverse planning problem in intensity-modulated radiation therapy (IMRT) [12–15] and the image reconstruction process in computed tomography (CT) [16,17]. In VAM and IMRT, the optimization attempts to create high dose regions in space for operation while minimizing side effects triggered by the undesirable dose elsewhere. In both cases, radiation emanates from an exterior source and passes through the subject. The non-negativity of radiation energy represents a common physical constraint of the design of radiation profile in these processes. For these reasons, some of the quantities considered in this text have direct counterparts in IMRT. Meanwhile, both the optimization problem in VAM and IMRT have benefited from the theoretical foundation established in computed tomography (CT). Work related to projection computation in VAM [1,18–20] references heavily the theory of CT and uses the analytical solution of CT to generate initial projection designs.

This work considers the projection optimization problem in typical VAM systems where material is converted by a single type of photoexcitation to deliver a single type of response. Extensions to scenarios with multiple types of photoexcitation [21] [22] are left for future work. The single-excitation-single-response systems under consideration generally share similar optimization goals, variables, and constraints. The general goal of the optimization is to seek a projection intensity distribution such that the material under exposure would achieve a spatial profile of conversion or a conversion-state-dependent property as close to a desired profile as possible. For conventional material systems [23,24], the desired conversion profile is often binary when only the geometry of the printed part is concerned—the local conversion either exceeds the gel conversion and the material polymerizes, or the dosage is insufficient and the resin remains unsolidified. In latest material developments, a novel class of "greyscale" photochemistries can exhibit drastically different material properties (such as elastic modulus and index of refraction) depending on the real (analog) degree of material conversion [25,26]. With this in mind, this work puts particular emphasis in improving analog control of conversion, which will enable VAM systems to leverage the above novel materials and rapidly produce functional multi-material devices. To describe the tomographic illumination accurately and improve property resolution, VAM systems need a generic light propagation model to account for optical effects such as attenuation, refraction, and scattering [24]. In all cases, the optimization model must also accommodate the physical constraints on variables, including non-negativity of physical irradiance and hardware limits such as maximum intensity and bit-depth. Although implementation varies across application contexts, this paper aims to identify the structures of the common goals and constraints such that users can systematically approach the projection optimization problem with classical mathematical tools.

Despite the common structure of the problem, the existing projection optimization schemes in VAM each constructed a different set of features to address different aspects of the problem. There lacks a consistent framework to bridge between these schemes and allow users to systematically understand, analyze and fine tune the parameters in the optimization. Table 1 tabulates the availability of various features in the three prominent optimization schemes considered in this study, namely Dose Matching (DM) [18], Penalty Minimization (PM) [18], and Object Space Model Optimization (OSMO) [19]. (Unless stated otherwise, all mentions of OSMO refer to its iterative scheme instead of its direct scheme designed for greyscale targets.) Among these schemes as reported, DM is the only scheme that iteratively optimizes a non-linear material response towards a real-valued profile. In contrast, the goal of both PM and OSMO is to fulfill a binary dose target by driving the dose values to meet a minimum or a maximum threshold. By not prescribing an exact dose profile to be met, this constraint satisfaction approach relaxes the problem and allows the dose values inside the target to be further separated from the dose values outside. In addition, the PM scheme further improves feasibility of the problem by forgoing control of the dose values at the edge of the binary target. In this paper, it will be shown that the explicit or implicit loss functions in these optimization schemes can be summarized into one general form which incorporates all the above features. Furthermore, the general loss function offers an additional parameter to control sparsity of the error distribution. This new optimization framework does not carry over the dose renormalization steps in the mentioned schemes. This design decision allows a consistent physical unit system to be maintained throughout the entire optimization (supplementary S.2).

**Table 1.** Loss function and availability of features in three prior optimization schemes and in the proposed scheme. Check mark ✓ and cross ✗ denotes presence and absence of features respectively. Remarks are made in-place.

|  | Dose matching (DM) [18] | Penalty minimization (PM) [18] | Object space model optimization (OSMO) [19]† | Proposed scheme |
|---|---|---|---|---|
| Loss function | $L_1$-norm | $L_1$-norm | No explicit objective function is stated in the original report. A $L_2$-norm loss function is discovered in this work | $L_p$-norm |
| Real-valued target | ✓ Real-valued | ✗ Binary | ✗ Binary | ✓ Real-valued |
| Non-linear material response | ✓Applied in loss function but not during initialization and metric evaluation | ✗ | ✗ | ✓Applied throughout initialization, optimization, and metric evaluation |
| Relaxation through soft constraint satisfaction approach | ✗ | ✓ | ✓ | ✓ |
| Relaxation by local weighting | ✗ | ✓ Edge de-emphasis | ✗ | ✓ |
| Global error sparsity control | ✗ | ✗ | ✗ | ✓ |

† Refers to OSMO iterative optimization scheme described in original report. Direct (non-iterative) approaches are generally considered as initialization steps in the current framework.

Currently there exist other projection computation schemes that also potentially fall within the same framework, but such schemes are outside the scope of main discussion. Bhattacharya *et al.* [18] showed that the heuristic optimization procedure used in the first reported tomographic VAM [1] is a special case of PM. Edge-Enhanced Object-Space Model Optimization (EE OSMO) [20] is an extension of OSMO designed to improve edge fidelity of binary targets by adding two additional model updates and reconstruction steps. Pazhamannil *et al.* also describe an iterative scheme [27,28] that uses similar steps such as initialization with filtered backprojection and gradient descent stepping as in other schemes.

Apart from optimization for binary targets, Rackson *et al.* [19] also describe a non-iterative approach to generate sinogram solution for real-valued targets. This direct approach is largely similar to the initialization procedures for the iterative OSMO for binary problems, with the exception of adding a constant offset after the frequency filtering step. Due to its relevance to initialization and its non-iterative nature, this approach is currently categorized as one of the possible ways to generate initial iterate.

Distinct from the above, Chen *et al.* [29,30] reported a projection optimization formulation that models the photon absorption as a stochastic process and maximizes the likelihood of a spatial profile of photon absorption. Although its optimization goal is still matching the delivered dose profile to the target profile, this formulation does not fit into the deterministic approach currently discussed here. Rigorous comparison between deterministic and stochastic models is left for further studies.

In summation, the projection optimization process in tomographic VAM would benefit from a more general and structured framework. Concretely, the value of the proposed loss function and optimization framework is to:

1. Provide deeper insights on existing optimization schemes by recasting and reinterpreting their parameters in the new formulation.
2. Enable a systematic and fair parameter study conducted through continuous single-parameter sweeps. Previous scheme-to-scheme comparisons are in fact switching multiple parameters at once.
3. Provide a large continuous parameter space where the trajectory of the optimization can be fine-tuned for specific applications. The proposed loss function generalizes existing features in previous schemes in a mutually compatible manner and provides additional material response modelling and error sparsity control features.

# 2 Methods

## 2.1 The generalized optimization model

In essence, the objective of the projection optimization problem is to locate a non-negative sinogram function $g$ that minimizes the error between the delivered dose response and the response target. The problem reads as $\min_{g \in S_{feasible}} \mathcal{L}$, where the loss function $\mathcal{L}$ quantifies the error as a soft constraint subjected to minimization and $S_{feasible}$ defines the set of solutions in sinogram function space $S$ satisfying the hard constraints. As will be discussed in supplementary section S.18, there are various types of hard constraints. For demonstration purposes, this paper provides examples where $S_{feasible}$ is the non-negative sinogram set $\{g \in S \mid g(\underline{r}') \geq 0 \; \forall \underline{r}'\}$. A key finding in this work is that a formulation of $\mathcal{L}$ is found to generalize the loss functions well in three seemingly disparate optimization schemes. The renormalization steps in these previously reported schemes are deliberately removed to preserve solution scale throughout optimization and allow hardware calibration with consistent physical units.

### 2.1.1 Forward model and loss function

In the optimization, the two crucial mathematical models that relate the sinogram $g$ to the dose response $f_m$ are the backpropagation model and the material response model. The backpropagation model represents the physical tomographic reconstruction process as a linear operation $P^*$ on sinogram $g(\underline{r}')$ in the computation of optical dose $f = P^*g$. A given backpropagation operator $P^*$ fully characterizes the (potentially shift-variant) optical impulse response $P^*\delta(\underline{r}')$ at each sinogram point $\underline{r}'$. The optical impulse response depends on the tomographic configuration and the material optical properties along the optical path of the corresponding beam. Taking the above optical dose $f$ as input, the material response model $\mathcal{M}$ simulates the dose response $f_m = \mathcal{M}(f)$ which is iteratively optimized towards the response target $f_T$. Table 2 provides further descriptions of these models. As a conceptual illustration, Fig. **1** graphically depicts the physical goal of optimization with example inputs and outputs.

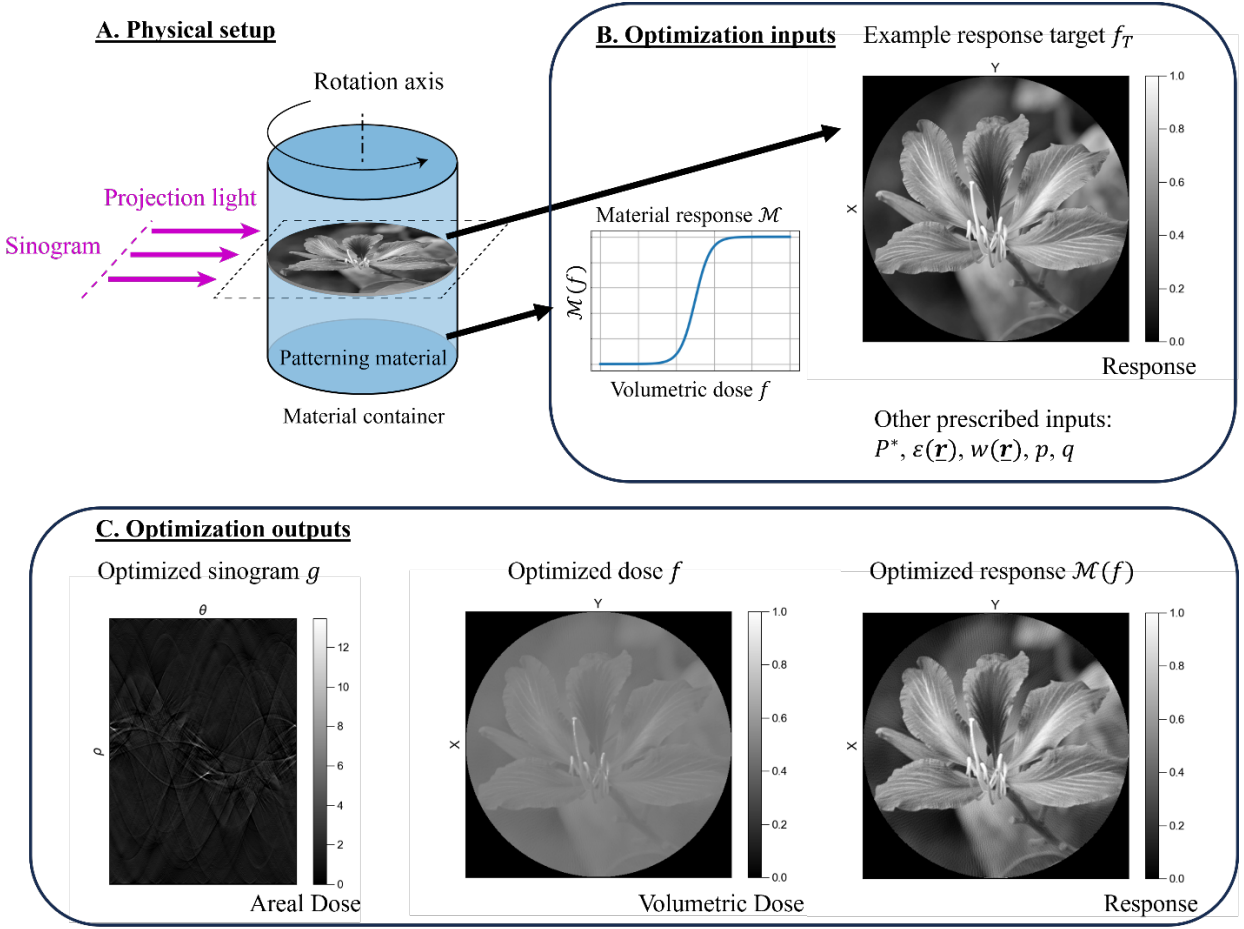

**Fig. 1.** Graphical depiction of the physical goal of the optimization problem in VAM. In the printing setup, the patterning material rotates relative to the projection light, as shown on the top left part of the figure (A). The user provides the response target, material response model, as well as other inputs listed top right part of the figure (B) to the optimization algorithm. The optimization iteratively updates the sinogram such that the tomographically reconstructed volumetric dose profile $f$ triggers a response profile $\mathcal{M}(f)$ as close to the response target $f_T$ as possible. A grayscale flower image is used as an example slice of an arbitrary real-valued response target with fine features. The photo of the flower is taken by Ianare [31] and is licensed under a CC BY-SA license. The bottom part of the figure (C) shows an optimized sinogram and its associated dose and response tomograms.

Taking an approach of soft constraint satisfaction, the loss function $\mathcal{L}$ only has a penalty term. Only when the local response value $\mathcal{M}(f)$ deviates from the response target value $f_T$ by more than response tolerance $\varepsilon$, does the loss function penalize the extent of the deviation $|\mathcal{M}(f(\underline{r})) - f_T(\underline{r})| - \varepsilon(\underline{r})$.

The new loss function takes the form of a weighted $L_p$-norm of the deviation evaluated over the soft-constraint-violation set $V = \{\underline{r} : |\mathcal{M}(f(\underline{r})) - f_T(\underline{r})| > \varepsilon(\underline{r})\}$ and raised to $q$-th power:

$$\mathcal{L} = \left( \int_V w(\mathbf{r}) \, \Big| |\mathcal{M}(f(\mathbf{r})) - f_T(\mathbf{r})| - \varepsilon(\mathbf{r}) \Big|^p \, d\mathbf{r} \right)^{\frac{q}{p}}. \tag{1}$$

The response target $f_T$, response tolerance $\varepsilon$ and weighting $w$ are all scalar fields parametrized with spatial coordinate $\mathbf{r}$. $p$ and $q$ are real scalars. Together with a known backpropagation operator $P^*$ and material response model $\mathcal{M}$, these parameters completely define the optimization problem. The physical meaning of these parameters is described in Table 2.

With the loss function defined, classical iterative optimization techniques can progressively lower the loss and find an optimal solution. The loss functions in each of the previously reported schemes falls under this general formulation $\mathcal{L}$. We name this generalized loss function the band constraint $L_p$-norm (BCLP) and the overall optimization problem as BCLP minimization. Fig. **2** graphically illustrates how the band constraint formulation generalizes both the real-valued and binary optimization targets in prior schemes.

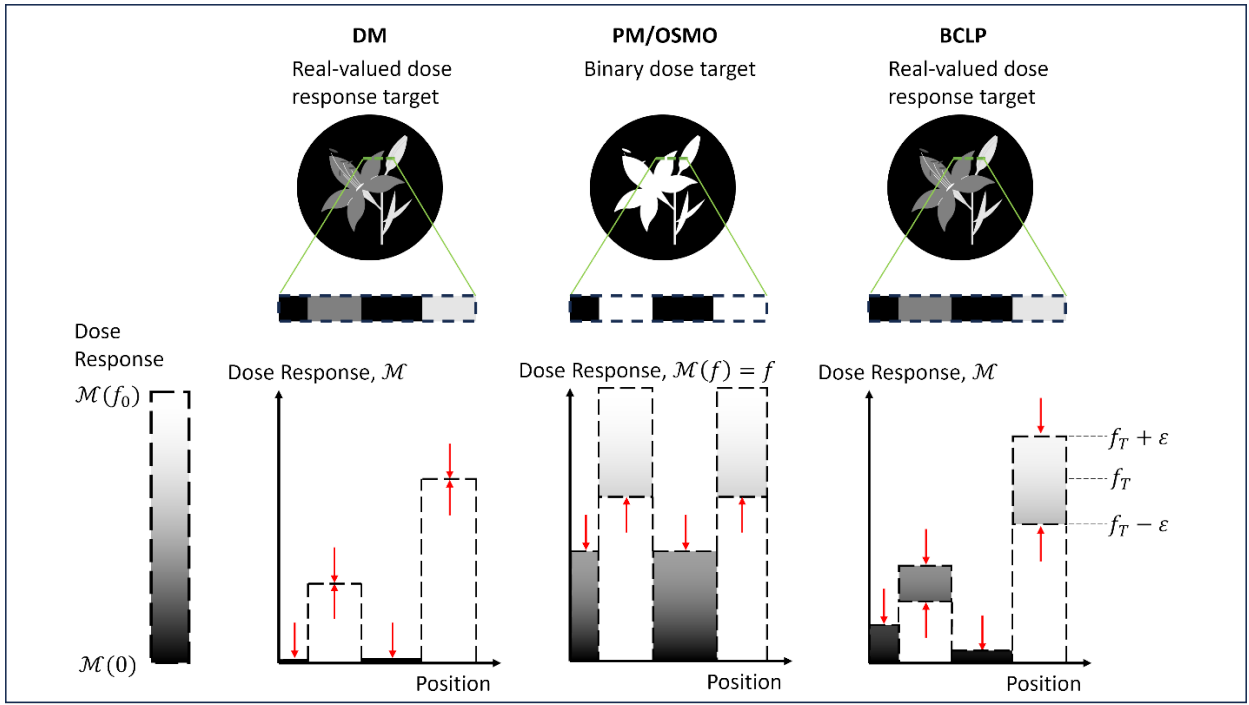

**Fig. 2.** Graphical illustration of the optimization goals in DM, PM, OSMO and BCLP. The red arrows indicate how the respective opitmization schemes drives the reconstructed dose or dose response profile towards a real-valued or binary target. The end of the arrows and the shaded regions indicate tolerance where the dose or dose response no longer penalizes the loss functions. In BCLP, the shaded tolerance bands are centered at $f_T(\mathbf{r})$ and with a interval of $2\varepsilon(\mathbf{r})$.

**Table 2.** Description and sources of inputs to the general optimization model.

| Input | Sources | Description |
|---|---|---|
| Backpropagation operator $P^*: S \to T$, where $S$ and $T$ are function spaces of square-integrable functions with compact support. | Assumed ideal propagation models (such as parallel beam and cone beam), optical propagation | $P^*$ is given and fixed during optimization. This linear operator maps each point on the sinogram to their 3D impulse response in the reconstructed tomographic dose. The linearity of $P^*$ in its argument directly follows from the theory of incoherent imaging. The shift-variance of the impulse response could arise from optical aberrations, |
| The domains of the sinogram functions in $S$ and tomogram functions in $T$ are three-dimensional spaces and are parametrized by sinogram coordinates $\underline{r}'$ and tomogram coordinates $\underline{r}$ respectively. | simulation (such as ray tracing and wave propagation in Fourier optics), or experimental characterization | non-parallel beam configuration, diffraction effects, attenuation, refraction, or local scattering.<br><br>The forward propagation operator $P$ is the adjoint operator of $P^*$. In matrix form, one is the conjugate transpose of another.<br><br>In IMRT, $P^*$ is referred to as the dose operator in a continuous model [32] and the dose deposition coefficient matrix in a discrete model [13,14]. |
| Material response, $\mathcal{M}: \mathbb{R} \rightarrow \mathbb{R}$ | Assumed or characterized material photochemical response | In this work, the material response model is an injective and differentiable function that maps optical dose $f$ to the dose response $f_m$ of the material. Examples of units of dose response are polymer degree-of-conversion, resultant elastic modulus and refractive index. Generally speaking, the temporal evolution of photochemical response can depend on the current material state and is influenced by transport phenomena. These further generalizations are not considered in this work. |
| Response target, $f_T(\underline{r}) \in \mathbb{R}$ | Prescription by user | The response target is a real-valued function of spatial coordinate $\underline{r}$. The optimization would attempt to steer the reconstructed dose $f(\underline{r})$ towards the tolerance band centered at $f_T(\underline{r})$. It takes the same unit as dose response. Examples in the current study use unity as the maximum $f_T$. |
| Response tolerance, $\varepsilon(\underline{r}) \in \mathbb{R}_{\geq 0}$ | Prescription by user | The response tolerance $\varepsilon(\underline{r})$ defines the width of the local tolerance band around the response target. The loss function does not penalize deviations within this band. A wide tolerance naturally relaxes problem constraints. $\varepsilon(\underline{r})$ shares the same units as $f_T(\underline{r})$ and $\mathcal{M}$. |
| Weights, $w(\underline{r}) \in [0,1] \subset \mathbb{R}_{\geq 0}$ | Prescription by user | The weights designate relative regional importance in the optimization. Non-uniform weights prioritize constraint satisfaction in the heavily weighted regions over lightly weighted regions. Users can adjust this emphasis according to their contextual priorities.<br>In order to strictly reproduce the alternate handling of overdosing and underdosing errors in OSMO, this weight needs to be defined differently in even and odd iterations. Supplementary S.12 and S.13 discuss this in detail. |
| $p \in \mathbb{R}_{>0}$ | Prescription by user | The $p$ value in the $L_p$-norm minimization controls the distribution of error (response deviation from the tolerance band). Minimizing $L_p$-norm at a particular $p$ value represents a trade-off selection (a point at Pareto front) between minimizing the spread of the non-zero error and minimizing the maximum error value. |

| | | A low $p$ value (0-1) emphasizes the occurrence of non-zero error and the corresponding minimization drives the error to a sparse distribution. |
| --- | --- | --- |
| | | A high $p$ value (1-∞) emphasizes the top percentiles of the error distribution and the corresponding minimization bounds the maximum error. |
| | | The special case of $p = 1$ puts equal emphasis on each region (voxel) of the spatial domain. The $p = 2$ case results in a Euclidean norm which would slightly emphasize high-error regions (voxels) in proportion to the error value. |
| | | In the discrete domain, functions are represented by tensors (or vectors, when flattened) and $L_p$-norm equals to $l_p$-norm. The special cases of $p \to 0$ and $p \to \infty$ are not investigated in this study. |
| $q \in \mathbb{R}_{>0}$ | Prescription by user | The value of $q$ changes the convergence behavior of the optimization. However, the locations of local and global minimizers on the loss function landscape are independent of $q$. This is because the ranking order of all solutions (with nonnegative loss) is preserved during exponentiation with positive power. |

Supplementary Table 3 list example choices of physical units for $g$, $P$, $f$, and $\mathcal{M}$.

### 2.1.2 Definitions and notations

In this text, the term "optimization formulations" refers specifically to the definition of loss functions which defines the optimization goal and dictates the location of local minima. The term "optimization schemes" refers to the high-level algorithmic procedures used in specific literature to arrive at a solution such as dose matching and penalty minimization. The term "optimization methods" is reserved for solution updating steps such as classical gradient descent and Newton's method.

To facilitate cross-referencing, notation in this work is chosen such that it is as close to prior work as possible while maintaining overall consistency.

Forward propagation operator $P$ and backpropagation operators $P^*$ are explicit generalization of forward projection and backprojection operators as defined in previous optimization schemes and in computed tomography. $P$ and $P^*$ are not necessarily modelling a parallel-beam tomographic configuration and may include modelling of scattering and refraction events. This naming is also to maintain proper distinction with the mathematical projection operation in projected gradient descent. They are not to be confused with forward propagation and backpropagation in context of optimization of neural networks. Supplementary S.3 and S.4 detail the composition of these propagation operators.

$S$ and $T$ in Table 2 refer to sinogram function space and tomogram function space, respectively. These infinite-dimensional function spaces are not to be confused with the three-dimensional (3D) domain of sinogram functions and tomogram functions. Sinogram coordinates $\underline{r}'$ and tomogram coordinate $\underline{r}$ parametrize the 3D domains of sinogram functions and tomogram functions, respectively.

Variables and operators in continuous form are not bolded. Tensors of order one and above are in bold. The number of bars under the tensor denotes the order of tensor. For example, there is one bar under vectors and two bars under 2D matrix.

Example parametrizations of tomogram function $f(\underline{r})$ and sinogram function $g(\underline{r}')$ are provided below.

Tomogram coordinate $\underline{r} = \begin{bmatrix} x \\ y \\ z \end{bmatrix}$, where $x$, $y$, and $z$ are Cartesian coordinates in 3D Euclidean space.

Sinogram coordinate $\underline{r}' = \begin{cases} [\rho \quad \theta \quad z']^T & \text{in parallel-beam configuration} \\ [\theta \quad \phi \quad z']^T & \text{in fan-beam configuration} \\ [\theta \quad \phi \quad \psi]^T & \text{in cone-beam configuration} \end{cases}$, where $\theta$ is the angular position of the projection gantry relative to the fixed material simulation domain. $\rho$ and $z'$ are Cartesian horizontal and vertical coordinates, respectively. $\phi$ and $\psi$ are the azimuthal and polar angles on the projection image relative to the optical axis of the projection gantry, respectively. The domain of sinogram functions can be constructed in cylindrical and spherical 3D Euclidean space in parallel-beam and fan-beam configuration respectively. The domain of sinogram function in cone-beam configuration can be constructed as a 3-sphere in 4D Euclidean space.

## 2.2 Initialization, analytical gradient, and iterative updates

While many derived first-order or second-order methods [33] can traverse the solution landscape and search for minima, this work uses projected gradient descent (PGD) [34,35] for demonstration. Compared to the simplest form of gradient descent, the PGD method additionally incorporates an operation in each iteration to project mathematically the current solution iterate onto the feasible set. The feasible set contains all solutions that satisfy hard constraints while the loss function penalizes the solutions that violate soft constraints proportionately.

### 2.2.1 Initialization

One common initialization approach is to project mathematically an optimum solution of the primal optimization problem without hard constraints $g_{uc,opt} = \arg\min \mathcal{L}$ to the feasible set $S_{feasible}$ to obtain an initial iterate solution $g_0$. In the general case where a direct solution of the unconstrained optima is not available, one must resort to iterative solution. Using iterative procedures to obtain such unconstrained optimization is undesirable for two reasons. Firstly, the unconstrained problem still needs its own initial guess solution. Secondly, optimizing the unconstrained problem would likely be as computationally intensive as optimizing the primal constrained problem because the former only omits projection onto $S_{feasible}$ which typically has negligible computational cost.

Instead of strictly solving $g_{uc,opt} = \arg\min \mathcal{L}$, we take $g_{uc,opt}$ as the solution of the simpler problem $\mathcal{M}(P^*g) = f_T$, where $g \in S$. This approach is motivated by the fact that any exact solutions $g_{exact}$ to $\mathcal{M}(P^*g) = f_T$, if they exist, would be precisely unconstrained optima as they would yield a minimum loss $\mathcal{L} = 0$. In particular, we solve for an approximate solution $g_{approx}$ to such an equation since an exact solution $g_{exact} = P^{*-1}\mathcal{M}^{-1}(f_T)$ may not exist. Solving the equation requires consideration of the right inversion of $\mathcal{M}$ and $P^*$.

The injective dose response function $\mathcal{M}$ introduced in this work is invertible over the range of $\mathcal{M}$ excluding the asymptotes. In practical settings, as long as the range of $\mathcal{M}$ covers most of the range of $f_T$, close approximation to $\mathcal{M}^{-1}(f_T)$ is readily available. Supplementary S.6 describes the analytical expression of $\mathcal{M}$, inverse $\mathcal{M}^{-1}$ and inversion of out-of-bound values. The illustration of optimization input in **Fig. 1** graphs an example dose response function.

Depending on the desired dose distribution $\mathcal{M}^{-1}(f_T)$ and the tomographic configuration represented by $P^*$, the equation $P^*g = \mathcal{M}^{-1}(f_T)$ may have many, one, or no exact solutions for $g$. Supplementary S.5 discusses some of the possible analytic and algebraic procedures to generate approximate solutions of $g$ for a given $\mathcal{M}^{-1}(f_T)$. Succinctly, the demonstration in this work approximates unconstrained optimum $g_{uc,opt}$ by:

$$g_{uc,opt} \approx g_{approx} = R_{Ram-Lak} P \alpha_{act,ab}^{-2} \mathcal{M}^{-1}(f_T), \qquad (2)$$

where $\alpha_{act,ab}(\underline{r})$ is the local absorption coefficients of the photoactive species responsible for the reaction and $R_{Ram-Lak}$ is the Ram-Lak filtering operation applied in the transverse coordinate of the sinogram. Then the initial sinogram iterate $g_0$ is computed by mathematically projecting the $g_{uc,opt}$ onto the feasible set:

$$g_0 = S_{feasible} g_{uc,opt} = \max(0, g_{uc,opt}), \qquad (3)$$

where $S_{feasible}$ is the projection operation onto the feasible solution set $S_{feasible} \subset S$. The feasible set $S_{feasible}$ is currently restricted to $\{g \in S \mid g(\underline{r}') \geq 0 \ \forall \underline{r}'\}$ due to non-negativity of optical areal power (intensity) and areal dose. In practical situations, additional hardware constraints on sinogram value such as quantization by digital bit-depth, upper bound by power limits or lower bounds by background level can be similarly applied in $S_{feasible}$.

### 2.2.2 Analytical gradient

The loss function has analytical gradient in the following form:

$$\nabla_g \mathcal{L}(\underline{r}') = q \mathcal{L}^{\frac{q-p}{q}} P\left(v(\underline{r})w(\underline{r}) E^{p-1} sgn\left(\mathcal{M}(f(\underline{r})) - f_T(\underline{r})\right) \frac{d\mathcal{M}}{df}\right)(\underline{r}'), \qquad (4)$$

where $v(\underline{r}) = \begin{cases} 1 & if \ \underline{r} \in V \\ 0 & if \ \underline{r} \notin V \end{cases}$ is the indicator function of $V$, $E = |\mathcal{M}(f) - f_T| - \varepsilon$ is abbreviation of response error, $sgn$ is the sign (signum) function and expression of $\frac{d\mathcal{M}}{df}$ depends on the choice of model $\mathcal{M}(f)$. This work uses a generalized logistic function as $\mathcal{M}$. Detailed derivation of the gradient $\nabla_g \mathcal{L}$ is provided in supplementary S.1. The expression of $\mathcal{M}$ and its derivative $\frac{d\mathcal{M}}{df}$ are provided in supplementary S.6.

### 2.2.3 Iterative updates

For ease of understanding and connection to existing optimization schemes, the demonstration in this study applies the projected gradient descent method. This optimization stepping method updates the iterate solution by

$$g_{k+1}(\underline{r}') = S_{feasible}\left(g_k(\underline{r}') - \eta[\nabla_g \mathcal{L}(\underline{r}')]_k\right), \qquad (5)$$

where $\eta$ is the step size and $S_{feasible}$ is the projection operation of the iterate onto the feasible set $S_{feasible}$. After $g_{k+1}(\underline{r}')$ is computed, the next loss function and loss gradient evaluation use $f_{k+1} =$ $P^*g_{k+1}$, and the loop continues. The optimization ends either if a maximum number of iterations is reached or if a convergence metric is met.

This work intentionally does not inherit the pre-iteration renormalization from OSMO. Removal of this normalization step preserves the scale of the tomogram and sinogram quantities throughout the iterative update process and allows the use of a physically meaningful unit system. Such scale preservation is strictly required when material response is not linear in optical dose. Supplementary S.2 discusses the benefits and necessity of a fixed physical unit. Supplementary S.9 provides example choices of the unit system. Section 2.3 discusses the quality metrics defined in this work. Supplementary S.16 demonstrates an optimization for a 3D binary response target.

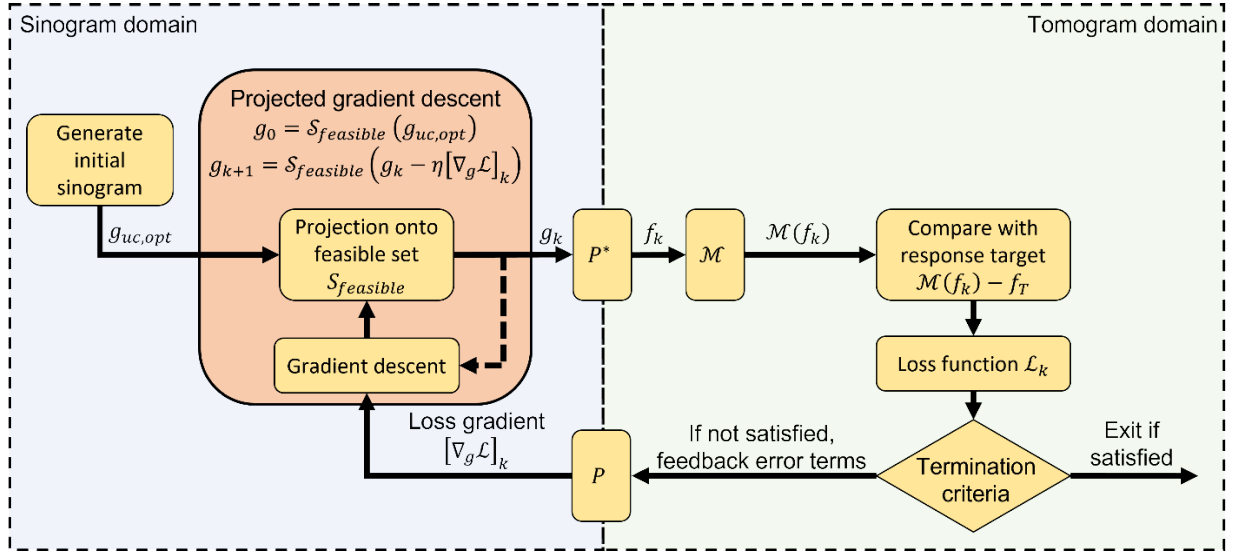

**Fig. 3.** Flow chart of projected gradient descent algorithm.

## 2.3 Performance metrics

### 2.3.1 Optimization metrics

This work uses only loss value as the main quality metric of solutions. Previous work [1,18,19] devised a variety of metrics such as Jaccard index, voxel error rate, and in-part dose range to represent favorable qualities of the solution. While they provide helpful information, these metrics are not necessarily optimized as they are not directly included in the loss function. In this work, we reference the loss value as the optimization metric and others as evaluation metrics. This is to signify that the latter are evaluated for informational purposes only and do not change the trajectory of the optimization. This distinction is practically important because it urges the users to align the formulation of the loss function to their desired solution characteristics. It is instrumental to establish mathematical correlations between these two types of metrics.

In practice, the loss function is not limited to modeling only one objective. In application contexts where competing goals exist, users can readily adopt a linear combination of multiple BCLP norms as the loss function ($\mathcal{L} = \sum_i \omega_i \mathcal{L}_i$, where $\omega_i$ is the weight for objective $\mathcal{L}_i$) such that these goals can be optimized in balance. Each of the norms $\mathcal{L}_i$ can evaluate the performance of the tomogram under different response functions (such as elastic modulus, degree-of-conversion, or refractive index), tolerance and regional emphases. Even more generally, other forms of loss functions can also be included. For the sake of simplicity, this multi-objective optimization is not included in the scope of this study.

### 2.3.2 Evaluation metrics

Despite the effort of designing intuitive loss functions, the loss value has limited descriptive power as it is only one scalar. Therefore, it is useful to generate additional function-specific information for performance evaluation and tuning of optimization parameters. For example, evaluation metrics are useful tools to observe the performance of certain spatial regions or the trade-off among multiple objectives in the loss function.

Previous work [1,18,19] has proposed a number of metrics to measure the accuracy and uniformity of dose reconstructions in binary printing. Currently these metrics cannot fulfill the additional measurement needs for reconstruction of real-valued targets. Supplementary S.8 discusses these additional needs and how various BCLP norms can perform such evaluation measurements.

# 3   Results and Discussion

## 3.1   Generalization of previous optimization schemes

The proposed formulation generalizes the optimization operations of the three previous optimization schemes to projection gradient descents on a band constraint loss function. The reformulation process provides new perspectives to understand and contrast the previous schemes. For conciseness, this section only highlights the key results of the reformulation (which is detailed in supplementary S.10, S.11, and S.12). To aid comparison, Table 3 tabulates the loss functions of the previous schemes, both in a simplified form and as a special case of BCLP. The variables in the loss functions are defined in the corresponding supplementary sections.

Table 3 The loss functions of previous optimization schemes in a simplified form and expressed as a special case of BCLP. The definition of the variables can be found in the corresponding supplementary sections (S.10 for DM, S.11 for PM, and S.12 for OSMO).

| Optimization scheme | Simplified loss function | Loss function expressed as a special case of BCLP |
|---|---|---|
| Dose matching (DM) | $\mathcal{L}_{DM} = \int \left| \sigma'(f(\underline{r})) - \Theta(\underline{r}) \right| d\underline{r}$ | $\mathcal{L}_{DM} = \int_V \left| \mathcal{M}(f(\underline{r})) - f_T(\underline{r}) \right| d\underline{r}$ <br> where $\mathcal{M}(f(\underline{r})) = \sigma'(f(\underline{r}))$ |
| Penalty minimization (PM) | $\mathcal{L}_{PM} = \rho_1 \int_{\sim V_1} (d_h - f(\underline{r})) d\underline{r}$ <br> $+ \rho_2 \int_{\sim V_2} (f(\underline{r}) - d_l) d\underline{r}$ | $\mathcal{L}_{PM} = \int_V w(\underline{r}) \left| |f(\underline{r}) - f_T(\underline{r})| - \varepsilon(\underline{r}) \right| d\underline{r}$ <br> where $w(\underline{r}) = \begin{cases} \rho_1, & \underline{r} \in R_1 \\ \rho_2, & \underline{r} \in R_2 \\ 0, & \underline{r} \notin (R_1 \cup R_2) \end{cases}$, <br> $f_T - \varepsilon = d_h \ \forall \underline{r} \in R_1$, and $f_T + \varepsilon = d_l \ \forall \underline{r} \in R_2$ |
| Object space model optimization (OSMO) | No explicit loss function | $\mathcal{L}_{OSMO,k} = \int_V w_k(\underline{r}) \left| |f(\underline{r}) - f_T(\underline{r})| - \varepsilon(\underline{r}) \right|^2 d\underline{r}$ |

|  |  | where $k$ is the iteration number, $$w_k(\underline{r}) = \begin{cases} 1 & if\ (k\ is\ even)\ and\ (\underline{r} \in OFP) \\ 1 & if\ (k\ is\ odd)\ and\ (\underline{r} \in IP), \\ 0 & otherwise \end{cases}$$ $f_T + \varepsilon = D_l\ \forall \underline{r} \in OFP$, and $f_T - \varepsilon = D_h\ \forall \underline{r} \in IP$ |
|---|---|---|

Among the three previous schemes, DM is the most straightforward to reformulate. It models the material response with a sigmoid function of dose and minimizes the integral of the absolute response error over the simulation volume. Its loss function is an unweighted $L_1$-norm. Since DM supports greyscale response targets, direct correspondence of the DM and BCLP variables can be made purely by factoring. Compared to DM, the BCLP formulation adopts a more general model of material response and applies the model consistently in optimization initialization and solution evaluation. Additionally, BCLP is equipped with local weighting and tolerancing for users to prioritize important regions and deprioritize inaccessible or less important regions. Also, the continuously adjustable $p$ value in the $L_p$-norm enables the user to control the sparsity of the error distribution.

The PM and OSMO iterative optimization take an approach different from DM. Their construction assumes a binary dose target from the outset. They attempt to keep the dose inside the part to be higher than a certain threshold and the dose outside of the target to be below another threshold. The PM and OSMO optimization proceed by correcting the violation of these single-sided (unilateral) soft constraints. In the reformulation, the dose thresholds in PM and OSMO are represented by one of the limits of the tolerance band in BCLP. By selecting the target value and the tolerance appropriately, the band constraint can effectively reproduce the behavior of the unilateral constraints in the prior schemes. This generalization unifies the correction approaches that are designed for greyscale and binary targets such that they can co-exist in the same optimization.

Under the general formulation, the constraint satisfaction approach of PM and OSMO can now be taken together with non-linear material responses and spatially variant weights. PM and OSMO minimize a $L_1$- and $L_2$-norm, respectively. Driven by their $p$ values, PM naturally favors a sparser error distribution than OSMO. In BCLP, these discrete choices of $p$ merges into a continuous variable that can be tuned to fit the task at hand.

OSMO uses an object space model as the optimization variable and algorithmically updates the model based on a sequence of operations that includes forward projection, truncation, backprojection, and renormalization. Hence, it does not have an explicit loss function in its original report. The reformulation of these algorithmic steps into projected gradient descent of the BCLP loss function is made possible by (1) rewriting the forward projections of the object space models as the sinogram variables, (2) leveraging the linearity of the projection (propagation) operation, and (3) identifying the equivalent model update steps in sinogram function space. Based on these observations, this work found a corresponding $L_2$-norm loss function which OSMO minimizes.

The reformulation made two subtle but valuable improvements to OSMO. Firstly, reformulating OSMO in projected gradient descent naturally maintains the scale of the iterating quantities and eliminates the need for the renormalization step in every OSMO iteration. Removing such renormalization steps and maintaining physical units are not only beneficial for setup calibration purpose but also crucial for handling non-linear material response. Supplementary S.2 discusses this particular aspect in greater detail. Secondly, the reformulation identified that the alternate handling of positive and negative errors in OSMO is optional. Although this alternate handling behavior can be reproduced in BCLP by defining weight $w_k$ as a function of iteration number $k$, a test performed in supplementary S.13 shows that the alternate handling is counterproductive. The optimization run with alternating handling converges slower than the one without alternate handling. Therefore, reformulating OSMO to BCLP offers a way to handle both positive and negative errors in parallel and potentially improve convergence.

## 3.2 Parameter study

The loss function completely dictates the location of local and global minima in the solution space. To obtain a desirable converged solution, it is important to choose values of loss function parameters such that they represent the application goal well. One of the benefits of the BCLP formulation is that it has a continuous parameter space which enables systematic parameter sweeps. In this section, a basic parameter study is performed to show the effect of each parameter and guide parameter selection. We study and analyze the effects of the following: the steepness of the material response $\mathcal{M}(f)$, the half-width of the tolerance band $\varepsilon(\underline{r})$ both globally and locally, the local weightings $w(\underline{r})$, and the $p$ value in the $L_p$-norm loss function. The parameter $q$ does not affect the location of global or local minima, and its effect on convergence properties is left for further study. Such convergence study should consider a multitude of factors including the choice of optimization method (such as gradient-based and quasi-Newton methods), initial solution, step size and convergence criterion.

The various parameters above are expected to have different influence on the local and global reconstruction accuracy. The goal of this study is to elucidate the effect of the above parameters and provide users with some intuition for fine-tuning the optimization result toward their contextual accuracy requirements. In practice, the reconstruction priority varies from application to application and strongly depends on the function of the part. For example, in an application of 4D printing[36,37] where the goal is to actuate a printed part by local differential swelling upon absorption of solvents, the variability of the degree-of-conversion limits the accuracy of the local swelling and the resulting motion. In this case, the accuracy requirements of the response (degree-of-conversion or swelling ratio) are dictated by the range of acceptable actuation distances. Furthermore, the user may want to prioritize accurate reconstruction of the folding mechanism over the other less critical regions. In many cases, the functional requirements of the part can be translated into a specification of response accuracy. Therefore, it is crucial for the user to set relevant accuracy priorities through the setting of loss function parameters.

In this study, the response target of each individual parameter sweep is chosen such that the effect of the parameter is obvious on plots. For conciseness, error histograms are only shown when they highlight the effect of the parameter. In each parameter sweep, only the sweeping variable is changed while all other settings are kept constant. As an exception, some optimization runs within the same sweep need to take a different step size to update meaningfully. These cases happen along the material response sweep and $p$ sweep. They are discussed in detail in the corresponding sections.

The convergence criterion for all parameter sweeps terminates the optimization when (1) the value of loss function reaches zero (a global minimum), or (2) the average-over-five absolute changes per iteration in loss function is less than 0.1% of the current value of loss function, as written in eq. (6) for iteration $k$.

$$\frac{\sum_{i=k-4}^{i=k}|\mathcal{L}_i - \mathcal{L}_{i-1}|}{5} \leq 0.001\,\mathcal{L}_k \tag{6}$$

Unless otherwise stated, the default parameter settings are: $p = 2$, $q = 1$, $\varepsilon = 0.05$, $w = 1$, and default material response settings as listed in Supplementary Table 1.

The response targets are 2D lying on a plane perpendicular to the rotation axis in tomographic scanning and they have 512 × 512 pixels. The lateral dimension $\rho$ of the corresponding sinogram has the same resolution as the lateral dimension of the response target. The sinograms have 360 tomographic projections spanning over 360 degrees, meaning the angular coordinate is discretized with 1-degree resolution. The sinogram values are represented by 32-bit floating point numbers, which have sufficiently high precision to limit bit-depth errors to the order of $10^{-6}$ units of response or below. The two real-valued response targets tested in the parameter sweeps are shown below in **Fig. 4**. The first response target is composed of four sinusoidal gratings with equal response amplitude but defined with different levels of quantization. The bit depths of the four gratings are increasing counterclockwise at 1, 2, 4, and 12 bits for the lower right, upper right, upper left and lower left grating respectively. The second response target is a grayscale picture of a flower with a blurry background. This target has complex structure down to pixel level resolution and has relatively even spread of grayscale values and spatial frequencies.

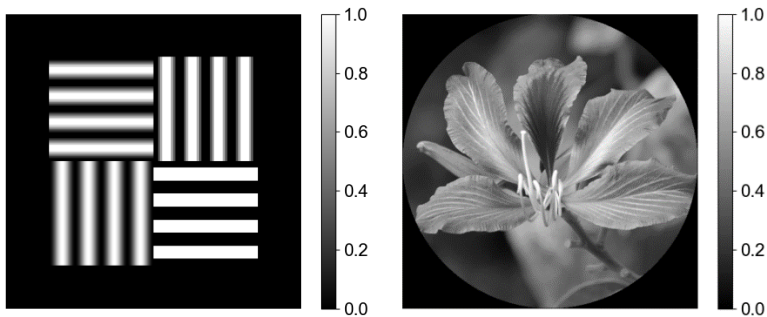

**Fig. 4.** The two response targets used in parameter sweeps. Left: Four sinusoidal gratings with 1, 2, 4, 12 bit-depths. Right: A flower. The photo of the flower is taken by Ianare [31] and is licensed under a CC BY-SA license.

### 3.2.1 Steepness of material response

Strictly speaking, the material response is not controlled by the user but instead governed by photochemical behavior of the material under consideration. Nonetheless, this parameter sweep aims to study the profound influence of the steepness of the material response on response reconstruction accuracy. This sweep sets the steepness parameter (loosely referred to as contrast) of the logistic material response function $\mathcal{M}(f)$ to be $B = [10, 25, 150]$ while keeping other parameters as default. To draw comparison with linear responses, an additional case is evaluated where $\mathcal{M}(f) = f$. All runs used a step size of 30 and converged, except the run $B = 150$ which used a different step size of 80 to produce meaningful results and did not converge under above criterion. This case is discussed at the end of this discussion. The material response curves, reconstructed dose responses, and dose response errors are plotted on **Fig. 5**. The presented solutions have loss function values of 97.6, and [53.5, 9.05, 49.3] for the linear case ($\mathcal{M}(f) = f$) and the three logistic response cases ($B = [10, 25, 150]$), respectively.

These results show that the contrast of the material is a strong determining factor for achievable reconstruction accuracy in response targets of various bit-depth. The overall response error is greatest in the run with a linear material response, and it is the smallest in the run with a logistic material response and $B = 25$. The better accuracy in the case with nonlinear logistic response is likely facilitated by the existence of saturation regions on both ends of response curve which allows the material to reject (or to be insensitive to) certain variances in dose. Comparatively, the run with $B = 10$ shows more error than the run with $B = 25$. This general observation is in line with the intuition in photolithography that a high-contrast material promotes response separation and rejects background exposure. Nevertheless, further increase in contrast from 25 to 150 actually leads to more error in regions with intermediate target response values. This provided a counterexample for the above intuition and suggests that a higher contrast is not always better. As expected, this near-binary response material is better suited to binary targets (having bit-depth of 1) and hence provided an almost perfect reconstruction for the grating on the lower right.

This result suggests that the optimal material response is target-dependent. If the material response model is parametric and differentiable (as implemented in this work), the parameters of the material response can be included in the optimization variables and be co-optimized towards specific target distributions under its respective constraints. As opposed to solely optimizing the sinogram with a preset $\mathcal{M}(f)$, this co-optimization could provide valuable information for experimental photoresist tuning and establish performance bounds for the chosen analytical form of response model.

The last run with logistic material response and $B = 150$ necessitated a step size different from other runs. The material response has drastically different response sensitivity at different dose values ($f$) and this large variance of local gradient ($d\mathcal{M}/df$) makes the reconstruction problem very ill-conditioned. With a relatively small step size, the optimization converges prematurely due to little change in the loss function. With a relatively large step size, the optimization lowers the loss stably at first but then the loss starts to oscillate and prevents convergence. We have not been able to locate a step size that both optimizes the response meaningfully and yet converges to the same criterion, and therefore present the results at 2000-th iteration. The convergence plot including all runs in this sweep is shown in supplementary S.15. To dampen such oscillation in practice, a gradient descent method with momentum [33] shall be considered.

The visible streaks on the reconstruction are likely aliasing artifacts that are generated during initialization by the Ram-Lak frequency filtering step. At the chosen angular discretization (1°/projection), the filtering step creates streaks that align with the gratings in the initial reconstruction. Although these artifacts remain visible in the converged solutions of the first two runs, they do not alter the conclusions of this study where the influence of material contrast dominates. Aliasing artifacts can be mitigated by using a tapering frequency filter instead of the Ram-Lak filter. Alternatively, they can be avoided completely by using an algebraic initialization method or a finer angular discretization that satisfies the Nyquist sampling criterion (0.448°/projection). This Nyquist sampling criterion is discussed in supplementary S.18.

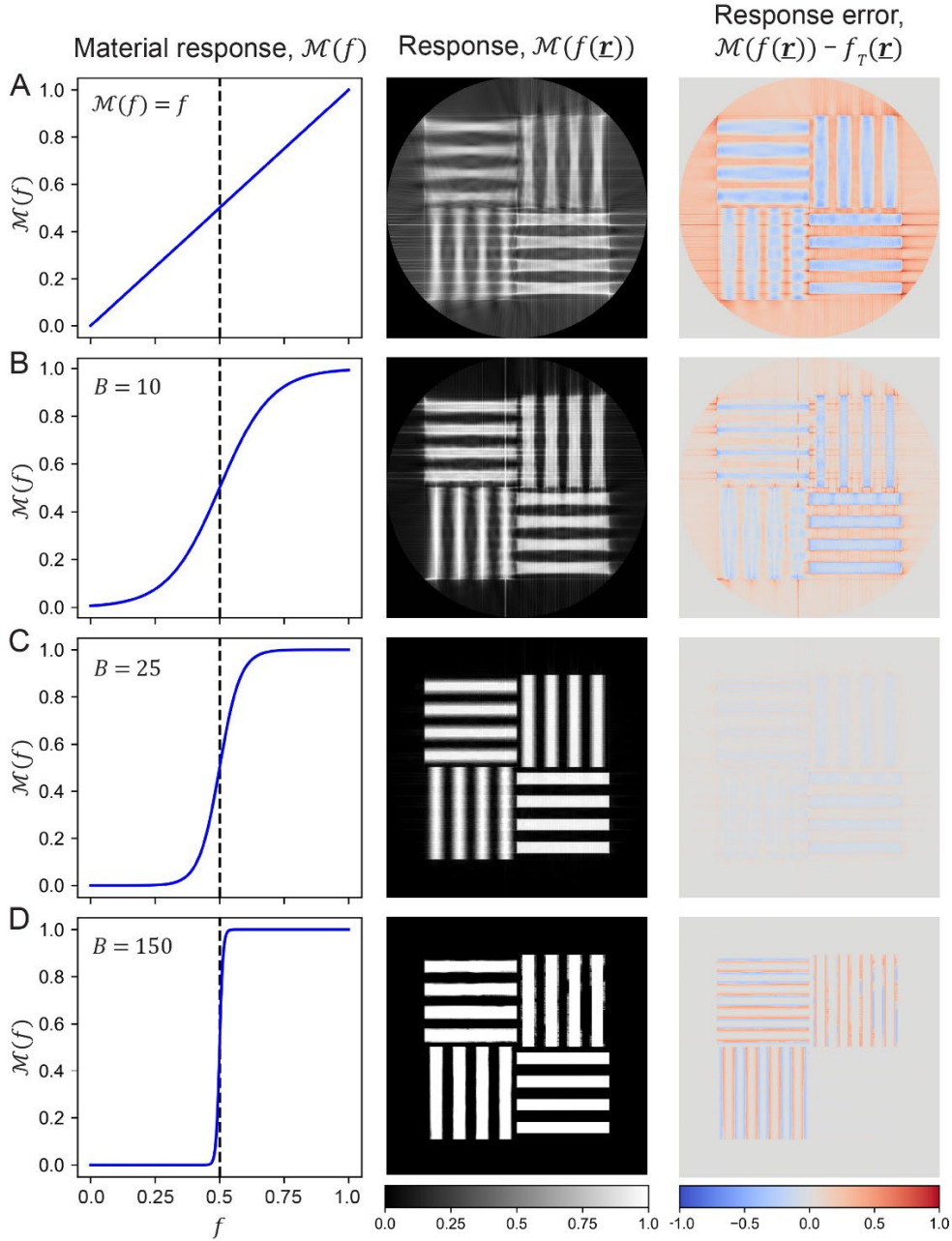

**Fig. 5.** Material response sweep. From top to bottom, the four rows (A-D) correspond to the four runs (linear, generalized logistic function with B=12,25,150) respectively. From left to right, the columns are plot of material response, response from tomographically reconstructed dose and response error from target respectively. The spatial axes are hidden for conciseness.

### 3.2.2  Global response tolerance

The soft constraint satisfaction approach relaxes the problem and directs the optimization effort to the regions where the response error exceeds the specified tolerance. The response tolerance specifies the cutoff of this non-linear correcting behavior and controls the amount of accepted error. This parameter sweep studies the effect of the global adjustment of the tolerance band width ($\varepsilon(\underline{r}) = \varepsilon_{global} \; \forall \; \underline{r}$) at four values of tolerances $\varepsilon_{global} = [0.2, 0.1, 0.05, 0.0]$. All optimization runs use default parameters (except tolerance), a step size of 10 and the flower image response target. **Fig. 6** shows the resulting response, response error and histogram of response error. The final loss values of these runs are [0, 0.879, 6.20, 18.9] respectively.

Apart from being a global minimum, a zero loss guarantees that the tolerance requirement is met everywhere and hence serves as an important stopping criterion in contexts with stringent performance requirements. Yet, it is not always physically possible to achieve a zero loss value. In the parameter sweep, the greater tolerance in the first run allows the optimization to achieve zero loss before convergence but it remains unclear whether the other three runs will ever achieve zero loss if they continue. As the tolerance tightens, it is fair to expect that the optimization has a lower chance of reaching zero loss and providing such a guarantee.

The rightmost column of **Fig. 6** shows that the response errors are distributed very differently among these cases. The first run with zero loss solution guaranteed that all voxels have a response error within the tolerance band (in grey shaded region) and there is a sharp peak in error population next to this limit (0.2 in absolute value). The formation of these peaks is a strong indication that the underlying error correction behavior stops right at the specified limit. The second run also shows similar peaks at the edge of the tolerance band but there are residue errors outside the band. In contrast, the third and fourth run shows no visible peaks near their respective 0.05 and 0 tolerance limits. Going from large to small tolerance runs, the population of voxels taking higher error grows. This trend agrees with the increasingly concentrated error on the response error plots in the middle column of **Fig. 6**. This phenomenon suggests that establishing an effective tolerance and accepting small errors allows the optimization algorithm to better allocate its effort on voxels with large errors. In other words, this form of problem relaxation enables the optimization to better focus on errors that are beyond the specified limits. Fundamentally, this observation helps to explain how PM and OSMO leverage this problem relaxation technique to deliver high contrast tomograms for binary printing.

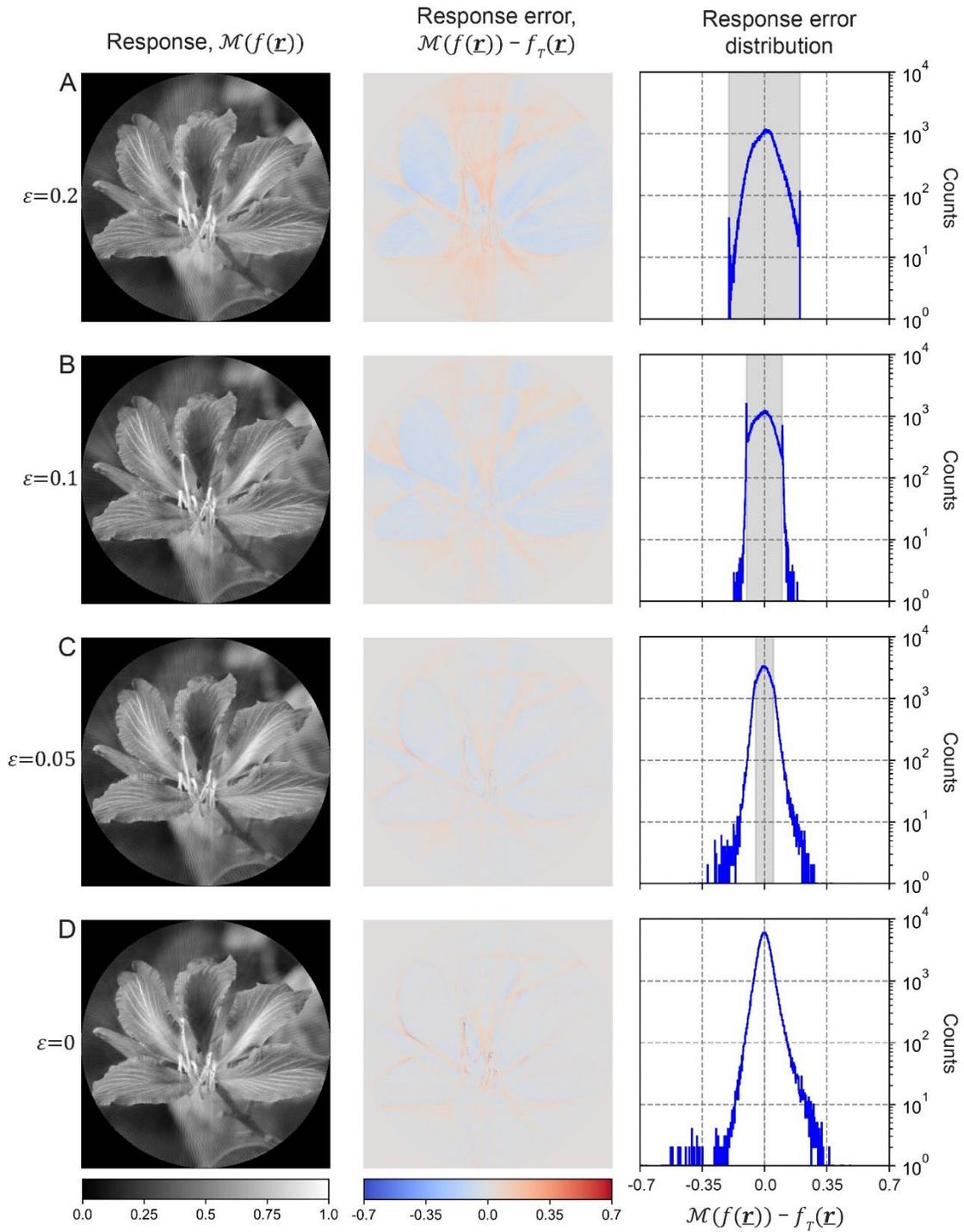

**Fig. 6.** Global tolerance sweep. From top to bottom, the four rows (A-D) correspond to the four runs with tolerance 0.2, 0.1, 0.05 and 0 respectively. From left to right, the columns are response resulting from tomographic reconstruction, response error, and histogram of error respectively. Gray shaded regions on the histogram of error indicate the width of the tolerance band. The spatial axes are hidden for conciseness.

### 3.2.3 Local tolerance and weighting

Both response tolerance $\varepsilon$ and weighting $w$ can be specified as a function of space to alter response accuracy locally. Nonetheless, they affect accuracy differently through different mechanisms and have different physical interpretations. The two following parameter sweeps aim to examine their individual function and discuss their differences. These two sweeps locally change the sweeping variable inside and outside a circular disk which has half of the diameter of that of the simulation domain. Both sweeps optimize for the "four grating" target and use default parameters unless specified. The step size in these sweeps is 100.

The local tolerance sweep applies the following tolerance value:

$$\varepsilon(\underline{r}) = \begin{cases} 0.4, 0.1, 0 & \underline{r} \in Disk \text{ for the three runs respectively} \\ 0.1 & \underline{r} \notin Disk \text{ for all runs} \end{cases} \quad (7).$$

**Fig. 7** plots the local tolerance, response, and response error for each of the three runs. Relative to the other runs, the first run leaves large errors inside the circle in the converged solution due to the much larger tolerance. The second and third runs produce almost identical results. This suggests that further refinement of tolerance beyond a certain value does not locally improve reconstruction accuracy.

The local weighting sweep applies the following weighting value:

$$w(\underline{r}) = \begin{cases} w_{disk} = \{0.1, 1, 3.5\} & \underline{r} \in Disk \\ w_{outside} = \{1.3, 1, 0.17\} & \underline{r} \notin Disk \end{cases} \quad (8)$$

for the three runs respectively. The weighting out of the disk is calculated by $w_{outside} = \frac{1 - 0.25(w_{disk})}{1 - 0.25}$ for the respective runs such that the spatial sum of the weightings applied in the simulation domain is constant across the three runs. This is intended to roughly maintain the scale of the problem relative to the fixed step size. **Fig. 8** plots the weighting distribution, response, and response error for each of the three runs. As intuitively expected, the regions with stronger weightings always have higher reconstruction accuracy relative to the regions with lighter weightings.

Despite their apparent similarity, setting the width of the local tolerance band is both mechanistically and practically different than setting the local weighting. Here we first point out the arithmetic differences. In the loss function, the relationship between the absolute response error and the tolerance is subtractive (in $V$). This subtraction locally shifts the optimization target for the absolute response error. Pictorially, the tolerance expands a soft equality constraint into a wider acceptance interval as depicted on Fig. **2**. In contrast, the relationship between the absolute response error and the weighting is multiplicative (after raising the former to the $p$-th power). One can also picture that the weighting term locally scales the differential volume $d\underline{r}$ (or a voxel in discretized form) in the loss function integral.

These arithmetic differences have practical implications. One obvious difference is that the weighting can almost infinitely scale the regional importance from zero to arbitrarily big numbers as an attempt to improve local reconstruction accuracy relative to other regions. In contrast, the tightest tolerance is zero and there is no tighter tolerance beyond zero to improve local reconstruction accuracy. Secondly, the tolerance and integration limit ($V$) together function as a truncated subtraction (or cut-off subtraction) which introduces more critical points in the loss function (with zero gradient). Once the voxel satisfies the tolerance, the contribution of that voxel to the loss gradient will vanish. This creates peaks in the error population away from zero error as shown in the first two rows of **Fig. 6**. In comparison, the effect of weighting is smooth and hence should not create peaks in error populations other than the zero-error bin. Therefore, weight and tolerance are not functionally interchangeable.

For the above reasons, there is practical value to have independent control over these two parameters. Heavily weighted regions do not have to pair with tight tolerance. As an example, heavily weighted regions with large tolerance can be used to establish hard limits that are easy to satisfy yet critical to part function. Conversely, light weights and small tolerance can be used to control the statistical average of less critical areas. These are scenarios that cannot be accommodated by just using one of the parameters.

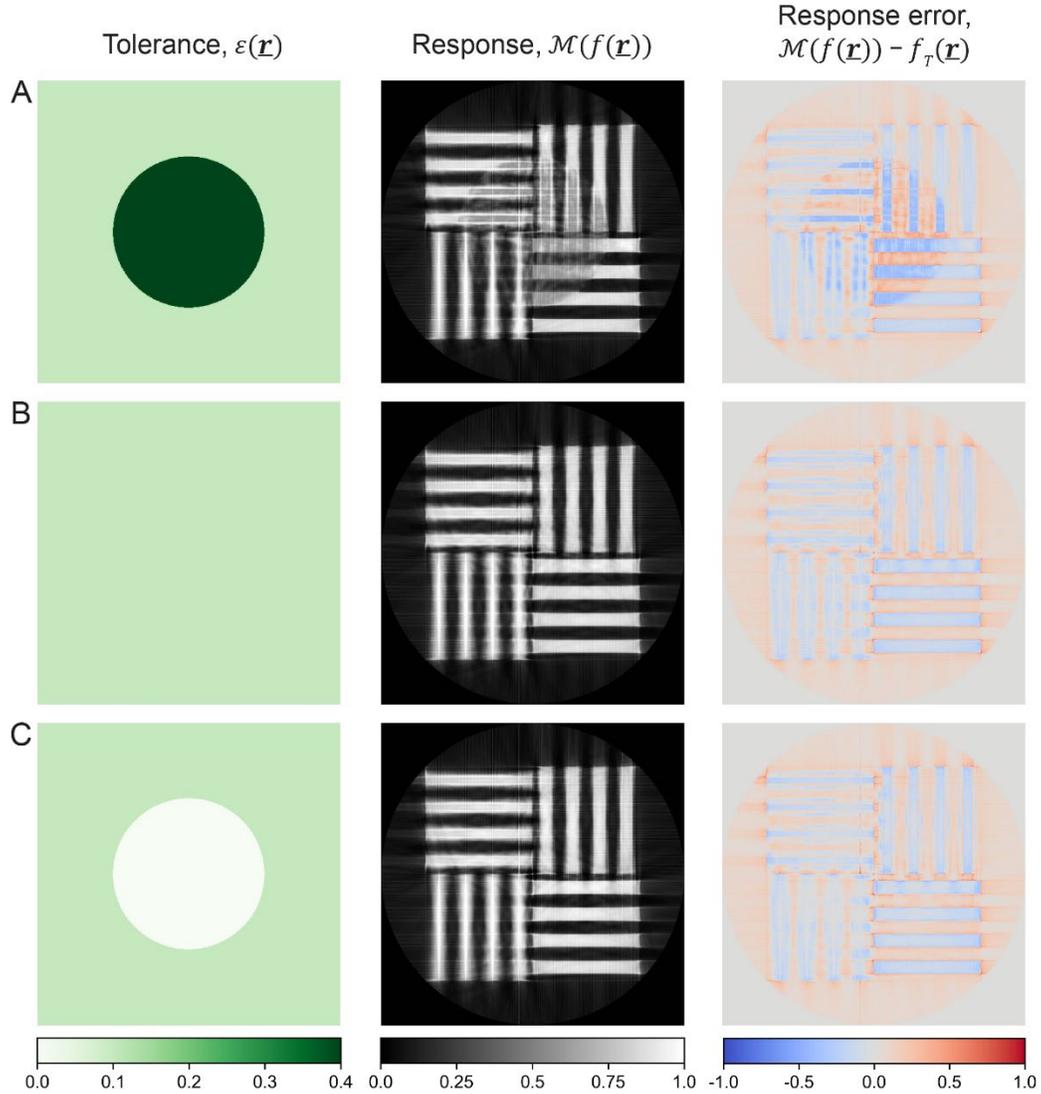

**Fig. 7.** Local tolerance sweep. From top to bottom, the three rows (A-C) correspond to three runs with different tolerance values $\{0.4, 0.1, 0\}$ inside a small disk. Tolerance outside the disk is kept at 0.1. The spatial axes are hidden for conciseness.

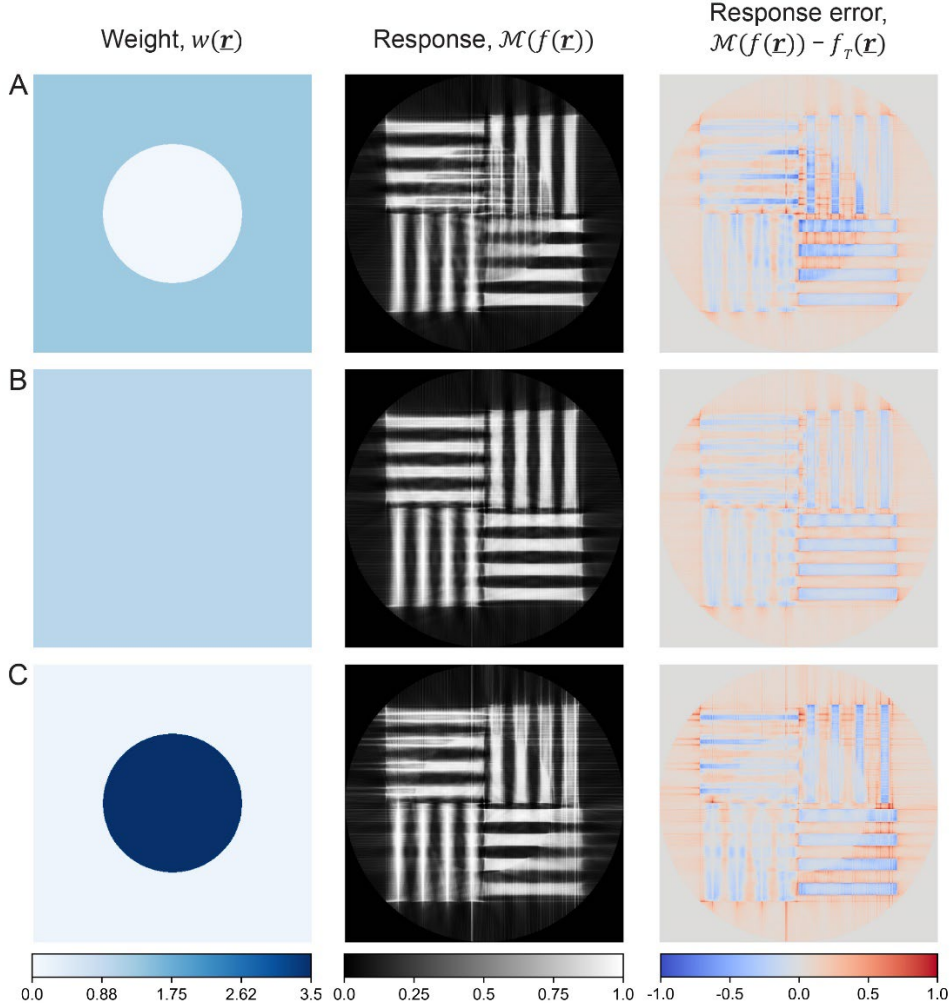

**Fig. 8.** Local weighting sweep. The three rows (A-C) correspond to the three runs with different weightings inside and outside of a disk. From top to bottom rows, the weightings inside are 0.1, 1 and 3.5 respectively and weightings outside are 1.3, 1, and 0.17 respectively. The spatial axes are hidden for conciseness.

### 3.2.4 Value of $p$

In the $L_p$-norm loss function, the $p$ value controls sensitivity of the loss function to the magnitude of the error $E(\mathbf{r}) = |\mathcal{M}(f(\mathbf{r})) - f_T(\mathbf{r})| - \varepsilon(\mathbf{r})$. A unity $p$ value leads to a loss function that depends linearly on the error at every point in space, while larger and smaller $p$ value biases the loss function towards the large and small error values, respectively. This parameter sweep shows the influence of the $p$ value on the resulting error distribution.

With other default parameters and the flower target, this sweep perform four optimization runs with $p = \{0.5, 1, 2, 20\}$. The four optimization runs use step sizes of $\{10^{-6}, 10^{-1}, 10, 10^2\}$, respectively, to accommodate magnitude difference of the loss function. **Fig. 9** shows the response, response error and histogram of response error of these runs.

From the histograms on the rightmost column of **Fig. 9**, the runs with smaller $p$ values show much higher population satisfying the tolerance (shaded in grey) than the runs with large $p$. This reflects that a small $p$ value prioritizes minimization of small errors (the lower percentiles) and helps to increase sparsity of soft constraint violation. On the other hand, the histograms also shows that the population with the largest errors diminishes as $p$ increases. This indicates that a large $p$ prioritizes the error minimization effort in large error regions (the upper percentiles). The response error plots in the middle column also confirm these findings as they show more concentrated errors with small $p$ and more spread-out errors with large $p$. Hence, tuning $p$ is an effective means to trade-off error sparsity with maximum range of error in practical situations.

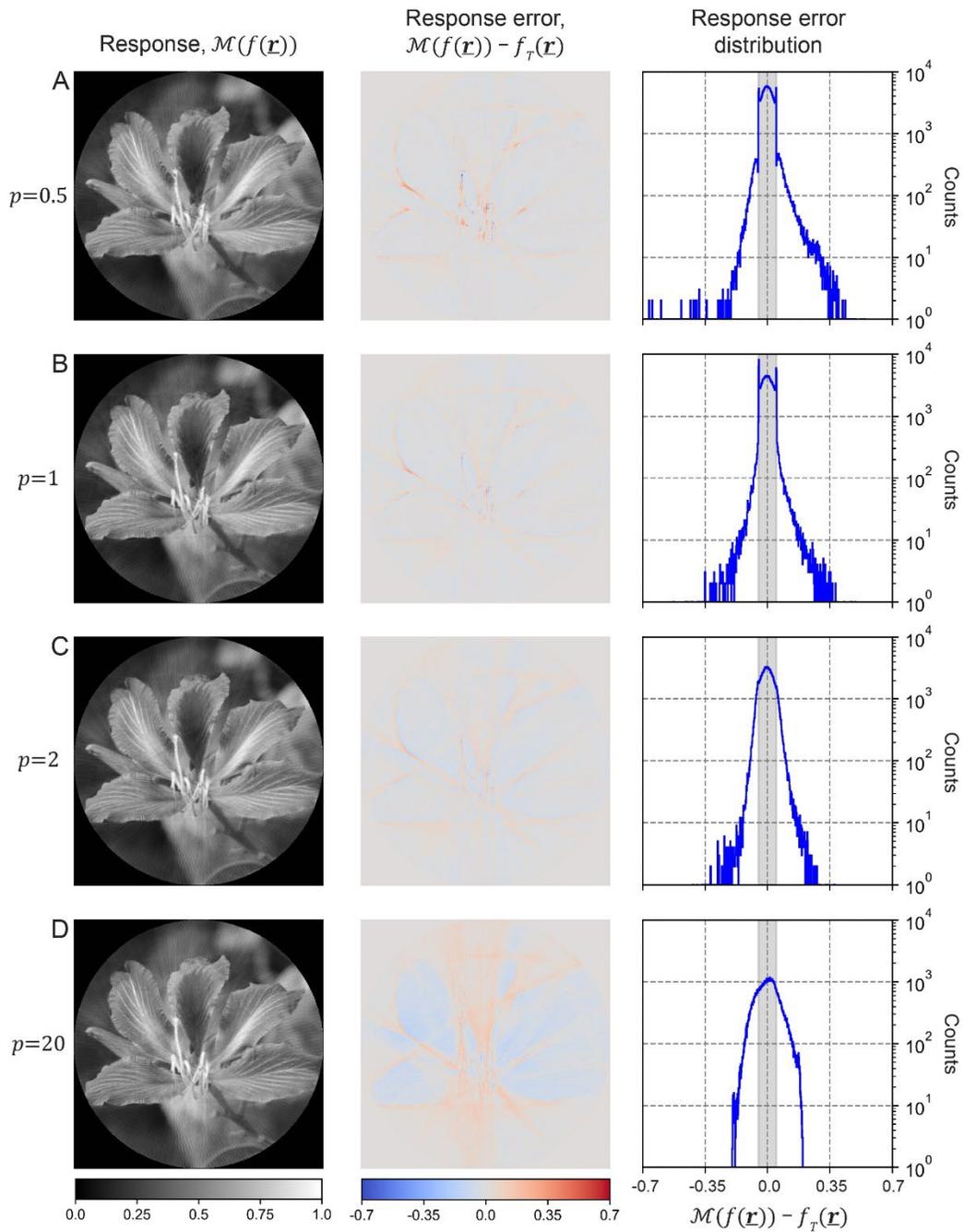

**Fig. 9.** $p$ sweep. The four rows (A-D) correspond to the four runs with $p$ value set to 0.5, 1, 2, 20 respectively from top to bottom. Gray shaded regions on the histogram of error indicate the width of the tolerance band. The spatial axes are hidden for conciseness.

# 4  Conclusion

This study uncovered the equivalent underlying operations of three recently published iterative projection optimization schemes and unified them in one generalized loss function. This reinterpretation offers fresh insights into existing optimization schemes and provides additional structures to consistently handle material response, preserve physical units, and control error sparsity. It presents a systematic approach to optimize tomographic projections for both binary and real-valued targets. The continuous and high-dimensional parameter space of this framework allows for the intuitive specification of manufacturing priorities for solution refinement under constrained settings. The results from the parameter study align with expectations and provide important application guidance. Promising directions for further generalization include incorporation of coherent propagation models, multi-wavelength material responses, and reaction–diffusion simulation.

# Code Availability

A Python implementation of the general projection optimization algorithm is available in the open-source repository https://github.com/facebookresearch/LDCT-VAM under GNU General Public License version 3.

# Tomographic projection optimization for volumetric additive manufacturing with general band constraint Lp-norm minimization

Chi Chung Li[1*], Joseph Toombs[1], Hayden K. Taylor[1], Thomas J. Wallin[2†]

[1]Department of Mechanical Engineering, University of California, Berkeley, Etcheverry Hall, 2521 Hearst Ave, Berkeley, CA 94720, USA

[2]Department of Materials Science and Engineering, Massachusetts Institute of Technology, 6-113, 77 Massachusetts Ave., Cambridge, MA 02139, USA

[*]Corresponding author: ccli@berkeley.edu

[†]Corresponding author: tjwallin@mit.edu

## S.1 Derivation and computation of analytic gradient

The generalized loss function is a weighted $L_p$ norm raised to $q$-th power. Let $E(\underline{r}) = |\mathcal{M}(f(\underline{r})) - f_T(\underline{r})| - \varepsilon(\underline{r})$. Acknowledging that $E(\underline{r}) > 0$ in $V = \{\underline{r} : |\mathcal{M}(f_r(\underline{r})) - f(\underline{r})| > \varepsilon(\underline{r})\}$ allows removal of the outer absolute sign of $|E(\underline{r})|$ and gives a simpler expression:

$$\mathcal{L} = \left( \int_V w(\underline{r}) \, E(\underline{r})^p \, d\underline{r} \right)^{\frac{q}{p}} . \qquad (SE.1)$$

With the chain rule of differentiation, the gradient of the loss function with respect to projection $g$ at a particular point $\underline{r_0}'$ in projection space is:

$$\frac{\partial \mathcal{L}}{\partial g(\underline{r_0}')} = q\mathcal{L}^{\frac{q-p}{q}} \int_V w(\underline{r}) \, E^{p-1} sgn\big(\mathcal{M}(f(\underline{r})) - f_T(\underline{r})\big) \frac{d\mathcal{M}}{df} P^* \delta^3(\underline{r}' - \underline{r_0}') \, d\underline{r} . \qquad (SE.2)$$

where $\delta^3(\underline{r}' - \underline{r_0}')$ is 3D Dirac delta in sinogram coordinates (which is equal to $\frac{1}{\rho}\delta(\rho - \rho_0)\delta(\theta - \theta_0)\delta(z - z_0)$ in cylindrical coordinates) and $sgn$ is the sign (signum) function. The quantity $P^*\delta^3(\underline{r}' - \underline{r_0}')$ is the impulse response (also called point spread function) of the backpropagation process at a sinogram point $\underline{r_0}'$.

Directly using the above $\frac{\partial \mathcal{L}}{\partial g(\underline{r_0}')}$ expression to compute the gradient would require an integral computation over real space for every point of interest $\underline{r_0}'$ in projection space. This approach is inherently computationally intensive. However, by using the property of the adjoint operator and the sifting property of the Dirac delta function, the above computation can be simplified such that the gradient with respect to all points in projection space can be evaluated in one forward propagation.

Rewriting the integral with the indicator function $v(\underline{r}) = \begin{cases} 1 & \text{if } \underline{r} \in V \\ 0 & \text{if } \underline{r} \notin V \end{cases}$ gives:

$$\frac{\partial \mathcal{L}}{\partial g(\underline{r_0}')} = q\mathcal{L}^{\frac{q-p}{q}} \int_\infty v(\underline{r}) w(\underline{r}) \, E^{p-1} sgn\big(\mathcal{M}(f(\underline{r})) - f_T(\underline{r})\big) \frac{d\mathcal{M}}{df} P^* \delta^3(\underline{r}' - \underline{r_0}') \, d\underline{r} . \qquad (SE.3)$$

The above expression is in fact an inner product in real space between $P^*\delta^3(\underline{r}' - \underline{r_0}')$ and the rest of the integrand $v(\underline{r})w(\underline{r})\, E^{p-1} sign\left(\mathcal{M}(f(\underline{r})) - f_T(\underline{r})\right)\frac{d\mathcal{M}}{df}$.

Since $P^*$ is the adjoint operator of $P$, then $\langle a, P^*b \rangle_T = \langle Pa, b \rangle_S$ for all functions $a$ and $b$ defined in tomogram function space $T$ and sinogram function space $S$ respectively. Then the gradient can be alternately expressed as an integral over the entire sinogram space:

$$\frac{\partial \mathcal{L}}{\partial g(\underline{r_0}')} = q\mathcal{L}^{\frac{q-p}{q}} \int_\infty P\left(v(\underline{r})w(\underline{r})\, E^{p-1} sgn\left(\mathcal{M}(f(\underline{r})) - f_T(\underline{r})\right)\frac{d\mathcal{M}}{df}\right)\delta^3(\underline{r}' - \underline{r_0}')\, d\underline{r}'. \quad (SE.4)$$

Then the sifting property of the Dirac delta function extracts the rest of the integrand at point $\underline{r_0}'$ and leads to:

$$\frac{\partial \mathcal{L}}{\partial g(\underline{r_0}')} = q\mathcal{L}^{\frac{q-p}{q}} P\left(v(\underline{r})w(\underline{r})\, E^{p-1} sgn\left(\mathcal{M}(f(\underline{r})) - f_T(\underline{r})\right)\frac{d\mathcal{M}}{df}\right)(\underline{r_0}'). \quad (SE.5)$$

The right hand side of the above equation is written in sinogram coordinates without integration. Writing the above expression in general sinogram coordinates $\underline{r}'$ gives:

$$\nabla_g \mathcal{L}(\underline{r}') = q\mathcal{L}^{\frac{q-p}{q}} P\left(v(\underline{r})w(\underline{r})\, E^{p-1} sgn\left(\mathcal{M}(f(\underline{r})) - f_T(\underline{r})\right)\frac{d\mathcal{M}}{df}\right)(\underline{r}'). \quad (SE.6)$$

Computationally, the optimization program first evaluates $v(\underline{r})w(\underline{r})\, E^{p-1} sgn\left(\mathcal{M}(f(\underline{r})) - f_T(\underline{r})\right)\frac{d\mathcal{M}}{df}$ in the tomogram domain and then directly forward propagates the evaluated quantity to obtain the gradient expressed in sinogram coordinates.

Overall, each iteration of the basic gradient descent algorithm involves a backpropagation to evaluate $f(\underline{r}) = P^*g(\underline{r}')$, a spatial integration to compute $\mathcal{L}$, and a forward propagation $P$ to perform a gradient update to $g$.

## S.2 Absence of renormalization per iteration

The OSMO scheme [1] and its variants [2] normalize the updated reconstruction dose $f_k$ by their maximum values in every iteration such that they do not exceed the numerical range of the response target, which is typically defined in the interval of $[0, 1]$. Consequently, the physical units of the variables are lost during the normalization step. In contrast, the BCLP formulation does not require such normalization to keep its response bounded. The BCLP loss function naturally penalizes the dose responses higher than the upper limit of the tolerance band and keeps them around the value of $f_T$. In fact, even when the BCLP tolerance band extends to infinity and functions as a unilateral soft constraint, the naturally saturating material response model would automatically bound the response values. Therefore, there is no reason for BCLP to inherit this normalization step from earlier approaches.

By eliminating such normalization, the proposed framework preserves scales and allows the entire optimization to be performed in physically meaningful units, which is important for optical and material calibration purposes. This scale-preservation feature is crucial in systems with nonlinear material responses $\mathcal{M}$, as the optical reconstruction dose $f_k$ must be consistently expressed in units accepted by the dose response model (for instance, $J/(cm^3)$). Otherwise, rescaling $f_k$ in every iteration would lead to inconsistent response evaluation. Therefore, it is necessary to remove this normalization step in general.

## S.3 Composition and implementation of propagation models in this work

In this work, the forward and backpropagation are abstracted as composite linear operators $P$ and $P^*$ respectively. As discussed in section S.9 about unit systems, these operators are constructed such that they directly map between volumetric quantities in real space (with units of $cm^{-3}$) and areal quantities in projection space (with units of $cm^{-2}$). The additional length unit comes from a multiplication with the absorption coefficient of the active species (which has units of $cm^{-1}$), denoted by $\alpha_{act,ab}(\underline{r})$.

In the weak attenuation regime, the linear attenuation coefficient of a mixture $\alpha_{total,at}$ is the sum of attenuation coefficients of all components $\alpha_{i,at}$. Each component attenuation coefficient $\alpha_{i,at}$ can further be broken down into the sum of the absorption $\alpha_{i,ab}$ and scattering coefficient $\alpha_{i,sc}$.

In the light transport model, both active and inert species could contribute to attenuation effects through absorption and scattering. However, only the light absorption by active species should count towards energy deposition that triggers dose response. Therefore, it is important to handle the active energy absorption term explicitly.

The volumetric dose is computed by multiplying the areal dose with absorption coefficient $\alpha_{act,ab}(\underline{r})$ in the tomogram domain. If we explicitly write out this multiplication in $P$ and $P^*$, we obtain $P = P_{at}\alpha_{act,ab}$ and $P^* = \alpha_{act,ab}P^*_{at}$, where $P_{at}$ and $P^*_{at}$ are respectively the forward and backpropagation operators that model all optical attenuation effects (computed with $\alpha_{total,at}$). $P_{at}$ and $P^*_{at}$ only map areal quantities to areal quantities (both with units of $cm^{-2}$). In the application of $P^*$ on some arbitrary sinogram distribution $g$, the multiplication of areal dose $P^*_{at}g$ with the local value of $\alpha_{act,ab}$ would give the excitation of the active species in volumetric dose units. One prime example of $P_{at}$ and $P^*_{at}$ is respectively the exponential Radon transform (ERT) and the associated exponential Radon transform (AERT) [3,4] referenced in prior optimization works such as PM and OSMO.

The demonstrations in this work use a custom ray-tracing implementation of the attenuated Radon transform (AtRT) [5,6] as $P_{at}$ and the adjoint of AtRT as $P^*_{at}$. Similar to ERT, AtRT models parallel-beam tomographic configuration in the regime of geometrical optics. Nevertheless, the AtRT generalizes the ERT and permits modeling of spatially varying attenuation coefficient $\alpha_{total,at}(\underline{r})$. Modeling attenuation coefficient as a function of space facilitates proper definition of the shape of the attenuating media (which is typically cylindrical) within the computation domain. More importantly, it also accommodates modeling of inhomogeneous media such as those with absorptive occlusions [1,7].

The demonstrations in this work compute $P$ and $P^*$ in a $512 * 512$ voxel grid at 500 voxel/cm sampling rate. In the square simulation domain, the total attenuation coefficient and active species' absorption coefficient only have support over the inscribed circular region:

$$\alpha_{total,at}(\underline{r}) = \alpha_{act,ab}(\underline{r}) = \begin{cases} 0.001 \; cm^{-1} \; if \; \underline{r} \in disk \; inscribed \; in \; the \; square \; simulation \; domain \\ 0 \; cm^{-1} \quad \quad \quad otherwise \end{cases}. (SE.7)$$

Apart from performing propagation directly, the implemented ray-tracing simulation can also produce an algebraic representation of this operation for faster subsequent evaluations. Except for the initialization steps described in section 2.2.1 and S.5.1, all propagation operations are performed using an algebraic representation of $P$ and $P^*$. The algebraic representation is discussed in section S.4.

## S.4 Composition, advantages and disadvantages of the algebraic representation of propagation operators

In simple tomographic configurations that can be approximated by the Radon transform or X-ray transform, the propagation operations can be evaluated relatively quickly. In contrast, optical propagation in complex tomographic configurations (such as those that involve scattering and refraction) may require costly optical simulations [8]. In these cases, it is more convenient to perform the backpropagation simulation only once to construct and store an algebraic representation of propagations $\underline{\underline{P}}^*: \mathbb{R}^{N_T \times N_S}$, where $N_T$ and $N_S$ are number of tomogram voxels and sinogram voxels respectively. In context of IMRT, $\underline{\underline{P}}^*$ is called the dose deposition coefficient matrix [9,10]. Conjugate transposition of the backpropagation matrix directly yields the corresponding forward propagation matrix $\underline{\underline{P}} = \left(\underline{\underline{P}}^*\right)^*$. These algebraic representations are matrices with entries being the coupling weights between the sinogram and tomogram voxels. Precisely, a particular column of $\underline{\underline{P}}^*$ represents the impulse response of a particular sinogram voxel in the voxelated tomogram domain.

Similar to the composition of $P$ and $P^*$ described in the last supplementary section (S.3), the matrices $\underline{\underline{P}}$ and $\underline{\underline{P}}^*$ are composed of

$$\underline{\underline{P}} = \underline{\underline{P}}_{at} \underline{\underline{\alpha}}_{act,ab}, \qquad (SE.8)$$

$$\underline{\underline{P}}^* = \underline{\underline{\alpha}}_{act,ab} \underline{\underline{P}}^*_{at}, \qquad (SE.9)$$

Where $\underline{\underline{\alpha}}_{act,ab}$ is a diagonal matrix with each diagonal element representing the local absorptivity of the active species. In discrete form, the forward propagation matrix operates on the discretized tomogram functions such as vectors $\underline{f}$ and $\underline{f}_T$. Correspondingly, the backpropagation matrix operates on the discretized sinogram functions such as vectors $\underline{g}_{uc,opt}$ and $\underline{g}_k$.

The algebraic representation $\underline{\underline{P}}^*$ of backpropagation has advantages over its simulation counterpart:

1. Direct storage of the algebraic presentations as matrices eliminates the need to repeat costly optical simulation during optimization. One backpropagation simulation can construct both matrices. Subsequent propagation operations can be carried out by efficient matrix–vector multiplications.
2. The matrix representation enables direct pseudoinverse or iterative least-square solutions to the problem $\underline{\underline{P}}^* \underline{g}_{uc,opt} = \mathcal{M}^{-1} \underline{f}_T$ in optimization initialization, which is discussed in the next supplementary section (S.5).

One of the limitations of the algebraic representation is that it requires a tremendous amount of memory even when it is stored in sparse format. Depending on the computational discretization of the tomogram domain and sinogram domain, even sparse matrices may not fit in the available random-access memory on personal computers.

Under a simplifying assumption that each sinogram voxel is coupled with a line through the tomogram domain with some non-zero weights, the number of non-zero elements in $\underline{\underline{P}}^*$ can be roughly approximated by the number of sinogram voxels $N_S$ multiplied by the number of tomogram voxels $N_{T,line}$ supporting a line in tomogram domain. For instance, a general 3D problem with $N_S = 512 * 512 * 180 \; voxels$ in sinogram space and $N_{T,line} = 512 \; voxels$ across an isotropic tomogram domain would require at least $24.2 * 10^9$ non-zero elements in $\underline{\underline{P}}^*$. If this $\underline{\underline{P}}^*$ is stored as sparse coordinate list (COO) format where its non-zero elements are stored as half-precision floating point (float16) and its matrix indices are stored as unsigned integer (uint32), it requires roughly 241.6 gigabytes of memory. For such a large problem, one can consider partitioning the matrix into blocks, storing them in secondary storage, and performing computation in sequence or over multiple cluster nodes.

3D problems that are shift-invariant along the rotation axis (such as parallel beam or fan beam configuration) can have propagation modelled in 2D slices normal to the rotation axis. If the propagation of the above problem can be completely represented by repeating a single slice, then $N_S = 512 * 180$ becomes much smaller and the memory requirement shrinks to 471.9 megabytes. Therefore, by leveraging the shift-invariant property in these problems, storage and handling of the algebraic representation $\underline{P}^*$ can be greatly simplified.

## S.5 Approximate solutions to the problem $P^*g = \mathcal{M}^{-1}(f_T)$ for initialization

At the start of the optimization, we project the approximate solution $g_{approx}$ of the equation $P^*g = \mathcal{M}^{-1}(f_T)$ onto the feasible set to obtain the first iterate solution $g_0$. This section discusses some of the possible analytical and algebraic methods to obtain such approximate solutions for a given distribution $\mathcal{M}^{-1}(f_T)$. Although both classes of methods are possible in the proposed optimization framework, all initial iterates in this work are generated with the analytical method described below.

The considerations in the selection of analytical and algebraic methods are highly similar to those in IMRT treatment planning. Users can take references from such related discussions [11].

### S.5.1 Analytical methods

Analytically, the exact solution can be written in terms of the right inverse of the backpropagation operator $P^{*-1}$ as $g_{exact} = P^{*-1}\mathcal{M}^{-1}(f_T)$. From the composition of $P^*$ described in section S.3,

$$P^{*-1} = \left(\alpha_{act,ab}P_{at}^*\right)^{-1} = P_{at}^{*\,-1}\alpha_{act,ab}^{-1}\ . \quad (SE.10)$$

This means that the operand $\mathcal{M}^{-1}(f_T)$ is first divided by the absorption coefficient of the active species $\alpha_{act,ab}(\underline{r})$ and then subjected to the inversion operation by $P_{at}^{*\,-1}$. With proper handling of $\alpha_{act,ab}(\underline{r}) = 0$ as discussed at the end of this section, the first process is computationally trivial. However, the inversion operation $P_{at}^{*\,-1}$ requires attention to the optical propagation model.

In this work, light propagation in the parallel-beam tomographic configuration is modelled as an attenuated radon transform (AtRT). Although there are several methods to invert AtRT [5,6,12–14], their computation is rather involved. Here, we approximate this inversion under the assumption of a weakly attenuating medium.

At the limit of negligible attenuation, the AtRT reduces to the conventional Radon transform and the right inversion of $P_{at}^*$ can be simply expressed in terms of a Fourier frequency filtering process followed by $P_{at}$ [15,1]. For a bounded tomogram function $f$,

$$P_{at}^* R_{Ram-Lak} P_{at} f \to f \quad as \quad \alpha_{total}(\underline{r}) \to 0 \ \forall \underline{r}, \quad (SE.11)$$

where $R_{Ram-Lak}$ represents a 1D Ram-Lak filtering operation applied on the Fourier frequency of the transverse sinogram coordinate. In other words, $R_{Ram-Lak}P_{at}$ is a right inverse of $P_{at}^*$ when attenuation is negligible.

Starting from the composition of $P^{*-1}$, and using the above approximated right inverse of $P_{at}^*$,

$$P^{*-1} = P_{at}^{*\ -1}\alpha_{act,ab}^{-1} \approx R_{Ram-Lak}P_{at}\alpha_{act,ab}^{-1} = R_{Ram-Lak}P\alpha_{act,ab}^{-2}. \qquad (SE.12)$$

Using this expression of $P^{*-1}$, we can construct $g_{approx}$ as:

$$g_{exact} = P^{*-1}\mathcal{M}^{-1}(f_T) \approx R_{Ram-Lak}P\alpha_{act,ab}^{-2}\mathcal{M}^{-1}(f_T) = g_{approx}. \qquad (SE.13)$$

In summary, if we apply this approximated solution $g_{approx}$ without enforcing any hard constraints, we should closely recover $f_T$ and get $\mathcal{L} \approx 0$.

$$\mathcal{M}(P^*g_{approx}) \approx \mathcal{M}\left(P^*P^{*-1}\mathcal{M}^{-1}(f_T)\right) = f_T \qquad (SE.14)$$

Taking this solution as $g_{uc,opt}$, we can start the projected gradient descent by projecting $g_{uc,opt}$ onto the feasible set $S_{feasible}$, compute loss and update solution with loss gradient.

It should be noted that the condition $\alpha_{act,ab} \neq 0$ for existence of $\alpha_{act,ab}^{-1}$ is always naturally satisfied inside the spatial region of interest $\{\boldsymbol{r}|w(\boldsymbol{r}) > 0\}$. In any physically meaningful photochemical or photothermal process, energy must be absorbed by the material to trigger a response. Therefore, $\alpha_{act,ab}$ must be non-zero and numerical infinities would not be generated inside the region of interest. In all physically-consistent settings where $\{\boldsymbol{r}|\alpha_{act,ab}(\boldsymbol{r}) = 0\} \cap \{\boldsymbol{r}|w(\boldsymbol{r}) > 0\} = \emptyset$, the region $\{\boldsymbol{r}|\alpha_{act,ab}(\boldsymbol{r}) = 0\}$ would not contribute to the loss function and therefore should not be included in the inversion. When numerically inverting $\alpha_{act,ab}$ in an array format, the numerical infinities resulting outside of the region of interest should be suppressed to be zero before the subsequent forward propagation operation. Alternatively, inversion of $\alpha_{act,ab}$ can be selectively performed at array entries corresponding to the region of interest.

### S.5.2 Algebraic methods

In the discrete form, the equation $P^*g = \mathcal{M}^{-1}(f_T)$ is written as a linear system of equations $\underline{\underline{P}}^*\underline{g} = \mathcal{M}^{-1}(\underline{f}_T)$ where $\underline{\underline{P}}^*$ is generally a rank-deficient non-square matrix, $\underline{g}$ is the vector form of $g$, and $\mathcal{M}^{-1}(\underline{f}_T)$ is the vector form of $\mathcal{M}^{-1}(f_T)$. Since $\underline{\underline{P}}^*$ naturally encapsulates $\underline{\underline{\alpha}}_{act,ab}$ and $\underline{\underline{P}}_{at}^*$ in one matrix, the algebraic solution process does not require consideration of the inversion of $\underline{\underline{\alpha}}_{act,ab}$ and $\underline{\underline{P}}_{at}^*$ individually. There are various algebraic methods to obtain a best approximation of $\mathcal{M}^{-1}(\underline{f}_T)$. This section gives an example of direct and iterative methods.

One of the direct methods is to apply the pseudoinverse (Moore–Penrose inverse) $\underline{\underline{P}}^{*\dagger}$ of matrix $\underline{\underline{P}}^*$ to vector $\mathcal{M}^{-1}(\underline{f}_T)$. This method naturally accommodates overdetermined and underdetermined linear systems. In the specific case where $\underline{\underline{P}}^*$ is square and full rank, $\underline{\underline{P}}^{*\dagger}$ coincides with the conventional square matrix inverse and $\underline{g} = \underline{\underline{P}}^{*\dagger}\mathcal{M}^{-1}(\underline{f}_T)$ gives the exact solution to the equation. For an overdetermined system, $\underline{\underline{P}}^{*\dagger}\mathcal{M}^{-1}(\underline{f}_T)$ gives the least-square solution to the problem $\min_{\underline{g}} \left\|\underline{\underline{P}}^*\underline{g} - \mathcal{M}^{-1}(\underline{f}_T)\right\|_2^2$, where $\|*\|_2$ is the Euclidean norm. For an underdetermined system, there are infinitely many solutions and $\underline{\underline{P}}^{*\dagger}\mathcal{M}^{-1}(\underline{f}_T)$ gives the minimum norm solution to the problem $\min_{\underline{g}}\|\underline{g}\|_2$ such that $\underline{\underline{P}}^*\underline{g} = \mathcal{M}^{-1}(\underline{f}_T)$. Singular value decomposition (SVD) has been used to compute the pseudoinverse and solve various kinds of tomography problems [16,17].

Since the matrix $\underline{\underline{P}}^*$ is typically very large, the direct method with SVD may be prohibitively expensive. In these cases, iterative methods with lower computational cost should be considered. Classically, well-studied algebraic algorithms in tomography such as Algebraic Reconstruction Techniques (ART), Simultaneous Iterative Reconstructive Technique (SIRT) and Simultaneous Algebraic Reconstruction Technique (SART) are designed to estimate the physical tomogram distribution from an experimentally measured sinogram distribution [18]. Though a swap of knowns and unknowns, these algorithms can equivalently be used to estimate the unknown sinogram $g_{approx}$ that would backpropagate into a known tomogram distribution $\mathcal{M}^{-1}(f_T)$. The advantages of these algorithms are that they only require low-cost matrix–vector products of the backpropagation matrix $\underline{P}^*$ but not its inversion or SVD which are much more computationally intensive. Particularly for least-square problems, algorithms such as LSQR [19] have been shown to provide computationally efficient and accurate results to tomography problems [20,21].

## S.6 Material response model

The material response model $\mathcal{M}$ captures the dose response $f_m$ to the optical reconstruction dose $f$. Using the preferred unit system as discussed in S.9, $f$ is the cumulative volumetric exposure defined as the time-integral of intensity multiplied by the absorption coefficient of the active species:

$$f(\underline{r}) = \int_{-\infty}^{\infty} \alpha_{act,ab}(\underline{r})\, I(\underline{r}, t)\, dt. \qquad (SE.15)$$

Generally, the material response (typically photopolymerization rate or material conversion rate) also depends on the intensity applied, diffusion of active species, and material state such as degree-of-conversion. However, accounting for these additional effects would require a reaction–diffusion simulation [22,23] which may not be easily differentiable. To limit the complexity of this demonstration, this work follows similar approaches taken by prior work and treats the response $f_m$ as a function of dose $f$ only. By parametrizing the response with dose, previous demonstrations successfully printed parts in a wide range of materials [24–29].

The model $\mathcal{M}$ used in this study takes the form of a generalized logistic function, commonly referred to as Richards's curve [30,31]. This chosen form is slightly more general than the logistic function (sigmoid) used in DM:

$$f_m = \mathcal{M}(f) = A + \frac{K - A}{\left(1 + e^{-B(f-M')}\right)^{\frac{1}{\nu}}}. \qquad (SE.16)$$

The effects of these parameters and the values used in current study are listed in Supplementary Table 1.

**Supplementary Table 1.** Parameters of the generalized logistic function.

| Parameter | Effect when other parameters are fixed | Default value in examples |
|---|---|---|
| $A$ | Left asymptote ($f \to -\infty$) | 0 |
| $K$ | Right asymptote ($f \to +\infty$) | 1 |
| $B$ | Steepness of the curve | 10 |
| $M'$ | Inflection point (value of $f$ that yields maximum slope) when $\nu = 1$. Generally, $M'$ shifts the curve left or right. | 0.5 |
| $\nu$ | Location of maximum slope relative to the two asymptotes. | 1 |

The corresponding first derivative is:

$$\frac{d\mathcal{M}}{df} = \frac{(K-A)\left(\frac{B}{\nu}\right)e^{-B(f-M')}}{\left(1+e^{-B(f-M')}\right)^{\frac{\nu+1}{\nu}}} = \left(\frac{1}{K-A}\right)^{\nu}\left(\frac{B}{\nu}\right)(\mathcal{M}-A)^{\nu+1}e^{-B(f-M')}. \qquad (SE.17)$$

In practice, it is possible to use other analytical functions or interpolation of purely numerical representations such as look-up tables. Nevertheless, strictly monotonic functions are numerically favorable because they are invertible and have non-zero gradient over the range of interest. Invertibility of $\mathcal{M}$ facilitates initialization of the optimization variable $g$. Functions with non-zero gradient everywhere would produce fewer saddle points on the loss function and are less likely to trap gradient update algorithms. The generalized logistic function has the above desirable properties and closely resembles the typical non-linear conversion response in photopolymerization processes, which often include induction and saturation periods near the beginning and end of polymerization respectively.

Earlier work (DM) only applies the material response model in optimization updates, but not in the initialization steps or evaluation metrics. The current work uses the following inverted response function eq. (SE.18) together with eq. (2) and eq. (3) to compute proper initialization $g_0$:

$$f = \mathcal{M}^{-1}(f_m) = M' - \frac{\ln\left(\left(\frac{K-A}{f_m-A}\right)^{\nu}-1\right)}{B}. \qquad (SE.18)$$

Section 2.3 and S.8 describe optimization and evaluation metrics that also consistently apply the response model.

This work implements the $\mathcal{M}^{-1}$ operation numerically as an interpolation. When it is queried at out-of-bound response values, the extrapolation dose values are taken to be the value at the bounds of the stored interpolant array.

## S.7 Nonconvexity of the loss function

The loss function $\mathcal{L}$ is generally non-convex with respect to the sinogram $g$ because the material response $\mathcal{M}$ is not a convex function. Therefore, local minima of $\mathcal{L}$ are not necessarily global minima.

However, there are special cases where $\mathcal{L}$ is convex. One of these cases is when material response is set as an affine function in dose (instead of being a non-convex logistic function), $p \geq 1$, $q \geq 1$ and $\varepsilon = 0$. In this case, the loss function is a convex $L_p$-norm of an error term $(\mathcal{M}(P^*g) - f_T(\underline{r}))$ that is affine in sinogram variable $g$. As it is known that a convex function of an affine function is overall convex, this special $\mathcal{L}$ is convex in $g$.

It should be noted that for the overall optimization problem to be convex, not only the loss function needs to be convex, the set of feasible solutions $S_{feasible}$ also needs to be a convex set. One such example of convex set is the set of non-negative sinogram $\{g \in S \mid g(\underline{r}') \geq 0 \; \forall \underline{r}'\}$ chosen in this work.

## S.8 Reconstruction evaluation metrics

Supplementary Table 2 lists a few evaluation metrics proposed in previous work for binary printing to measure volume of segmentation error and dose uniformity. While these existing metrics are good for their intended purposes, they cannot fulfill several additional needs in evaluating reconstructions for real-valued response targets. These additional needs are:

1) volume segmentation based on a spatially variant constraint,
2) segmentation and error evaluation in material response units, and
3) statistical measures for magnitude of response error.

To address these needs, this section discusses a few selected forms of the BCLP norm (independent of the loss function) that can provide reconstruction quality information.

**Supplementary Table 2.** Some reconstruction evaluation metrics in previous work.

| Metric | Description | Mathematical expression |
|---|---|---|
| Jaccard Index ($JI$)[26,32] | A spatial similarity metric, evaluated as the intersection over union of the binarized dose and the binary target | $JI = \frac{V_B \cap V_T}{V_B \cup V_T}$, where $V_B$ is the binarized response set and $V_T$ is the binary target set. |
| Voxel Error Rate ($VER$) [1] | An error metric, evaluated as the normalized overlap between the histograms for in-part (IP) and out-of-part (OFP) region. | $VER = \frac{W}{N}$, where $W$ is the number of out-of-part voxels that exceeded the lowest dose level of in-part voxels, and $N$ is the total number of voxels in the simulation volume. |
| In-Part Dose Range ($IPDR$) [1] | A uniformity metric, evaluated as the range of dose variation inside the IP region. | $1 - d_{IP,lowest}$ Where $d_{IP,lowest}$ is the lowest dose level in IP voxels in normalized dose units. |

Being designed for binary printing, the previous evaluation metrics segment the volume through a binarization operation. This binarization aims to represent the physical solid-liquid phase separation of the material during the development step in part post-processing. The binarization process segments the volume into a gelled population and an ungelled population by applying a spatially invariant dose threshold. Such dose binarization thresholds are determined by Otsu's method (for JI) [26,33], by maximizing the Jaccard index [32], or taken as the minimum IP dose (for VER and IPDR) [1]. In contrast, reconstruction of a real-valued response may require separation of more than two response levels. For example, certain mechanical metamaterials may call for five different levels of elastic modulus to be well distinguished between each other. Or, optical elements with gradients of refractive index may need the printed part to reproduce a continuously varying index profile. In general, applications may need a metric to represent how well a real-valued response reconstruction follows a multi-level or infinite-level design.

Furthermore, functionally relevant metrics should also segment volume and measure error in units of response instead of dose, since the response–dose relationship is not necessarily linear. As reported, the previous metrics do not have a dose–response error mapping built-in. A large error in dose may be a small error in response, and vice versa. Although this distinction between dose and response fades away when binarizing strictly monotonic responses for estimating phase separation, measuring response error is particularly important in the more general context with real-valued response targets. Finally, the existing metrics pay little attention to the magnitude of response error but mostly emphasize either the volume error of segmented regions (as in JI and VER) or the maximum error in a particular segmented region (as in IPDR).

To address the aforementioned needs, users can formulate custom BCLP norms as evaluation metrics independently from the loss function. The BCLP norm naturally performs volume segmentation with a spatially variant and infinite-level tolerance band. The band constraint classifies spatial regions to be either out of or within tolerance (in $V$ or its complement). Users can prescribe regions of interest by setting $w(\underline{r})$ as an indicator function. All metrics formulated in BCLP naturally segment volume and evaluate error using the material response $\mathcal{M}$ which is consistent with the loss function. The value of $p$ controls the sensitivity of the metric towards the magnitude of response error.

In the following examples, we define the shorthand:

$$\mathcal{L}_{\substack{p \to p_0 \\ q=q_0}} = \lim_{p \to p_0} \mathcal{L}|_{q=q_0} . \quad (SE.19)$$

For instance, $\mathcal{L}_{\substack{p \to 0^+ \\ q=p}} = \int_V w(\underline{r}) \, d\underline{r}$ gives the constraint violation volume over the region selected by indicator function $w(\underline{r}) \in \{0,1\}$. Analogous to JI and VER, $\mathcal{L}_{\substack{p \to p_0 \\ q=q_0}}$ provides volume information of the constraint violation region but it is defined over a multi-level or infinite-level tolerance band. $\mathcal{L}_{\substack{p \to p_0 \\ q=q_0}}$ is more closely related to VER than to JI because both $\mathcal{L}_{\substack{p \to p_0 \\ q=q_0}}$ and VER are proportional to the volume that violates a certain constraint.

On the other hand, $\mathcal{L}_{\substack{p \to \infty \\ q=1}} = \max_{V \cap \{\underline{r}:w(\underline{r})=1\}} \left( \left| \mathcal{M}(f(\underline{r})) - f_T(\underline{r}) \right| - \varepsilon(\underline{r}) \right)$ gives the maximum absolute error away from the tolerance limit in the region where the constraint is violated and indicator function $w(\underline{r}) \in \{0,1\}$ equals 1. A zero tolerance $\varepsilon(\underline{r})$ will leads to a metric that measures maximum absolute error from the $f_T(\underline{r})$. This region-specific metric can provide information on the range of response error similar to the maximum dose error in IPDR. Formally, the infinity-norm should be expressed in essential supremum instead of a maximum [34], but this technicality is not of concern in practical scenarios where computation is performed in discrete domains.

In addition, $\mathcal{L}_{\substack{p \to 1 \\ q=1}}$ and $\mathcal{L}_{\substack{p \to 2 \\ q=1}}$ gives the Manhattan norm and Euclidean norm of the absolute response error, respectively, measured from the tolerance limit, weighted by $w(\underline{r})$ and evaluated over the constraint violation region $V$. Roughly speaking, $\mathcal{L}_{\substack{p \to 1 \\ q=1}}$ is proportional to the weighted mean of the absolute response error and $\mathcal{L}_{\substack{p \to 2 \\ q=1}}$ is proportional to the weighted root-mean-square response error. As mentioned earlier, the $p$ value controls the sensitivity of the metric towards the magnitude of the response error. $\mathcal{L}_{\substack{p \to 1 \\ q=1}}$ and $\mathcal{L}_{\substack{p \to 2 \\ q=1}}$ have linear and quadratic contributions, respectively, from the response errors in $V$.

Finally, the proposed BCLP evaluation metrics with $0 < p < \infty$ are differentiable and can be readily included in the loss function as additional objectives for optimization. Users can use $p \ll 1$ and $p \gg 1$ to approximately optimize for the excluded edge cases of $p = 0$ and $p \to \infty$, respectively.

# S.9 Possible choice of physical units and difference between volumetric and areal dose

Example combinations of physical units are tabulated non-exhaustively below. The user-defined material response can be degree-of-conversion, $DOC$ [$unitless$] or its correlated material properties such as stiffness, $E$ [$Pa$] and refractive index, $n$ [$unitless$]. In the following table, we simply leave the unit of material response unspecified. This work adopts the unit combination in the first row of Supplementary Table 3.

**Supplementary Table 3.** Possible choice of physical unit systems.

| $g$ | $P$ and $P^*$ | $f = P^*g$ | $\mathcal{M}, f_T$, and $\varepsilon$ |
|---|---|---|---|
| areal dose [$\frac{J}{cm^2}$] | Absorption coefficient [$\frac{1}{cm}$] | Cumulative volumetric dose [$\frac{J}{cm^3}$] | material response to volumetric dose |
| intensity [$\frac{W}{cm^2}$] | Time*absorption coefficient [$\frac{s}{cm}$] | Cumulative volumetric dose [$\frac{J}{cm^3}$] | material response to volumetric dose |
| areal dose [$\frac{J}{cm^2}$] | None (Defining $P = P_{at}$, and $P^* = P^*_{at}$ using notations in section S.3) | Cumulative areal dose [$\frac{J}{cm^2}$] | material response to areal dose |
| intensity [$\frac{W}{cm^2}$] | Time [$s$] (Defining $P = P_{at}$, and $P^* = P^*_{at}$ using notations in section S.3) | Cumulative areal dose [$\frac{J}{cm^2}$] | material response to areal dose |

Relative to the unit system in the first row, the system in the second row moves the time dimension into propagation operators and serves as a more natural choice of units when tomographic scanning happens at non-uniform scan rates. The third and fourth rows only differ from the first two by omitting the multiplication with absorption coefficient $\alpha_{act,ab}(\underline{r})$ in propagation operators $P$ and $P^*$.

While all four combinations are programmatically possible, the usage of unit systems in third and fourth rows are discouraged for two reasons.

Firstly, using reconstruction quantity $f$ in unit of cumulative areal dose ($[\frac{J}{cm^2}]$) is much less relevant than cumulative volumetric dose ($[\frac{J}{cm^3}]$) (also referred to as absorbed optical dose[25]) in the context of photochemical reactions. The first law of photochemistry (Grothus-Draper Law) states that light must be absorbed by a chemical substance for a photochemical reaction to take place. It is the absorbed portion of light being responsible for the reaction. Therefore, the unabsorbed portion of light in cumulative areal dose is photochemically irrelevant. Despite its popularity in literature, the equipment-oriented areal dose is less relevant than the volumetric dose when considering the actual photoexcitation and material calibration.

Secondly, the unit combinations in the third and fourth rows do not explicitly consider the absorption coefficient of the active species $\alpha_{act,ab}$ and therefore cannot accommodate the scenarios where this coefficient varies spatially. The simulation volume does not necessarily have the active material uniformly distributed everywhere. The need for a spatial description of $\alpha_{act,ab}$ is very similar to the need for a spatial description of total attenuation coefficient $\alpha_{total,at}$ as discussed in supplementary S.3.

## S.10 Generalization and interpretation of the dose matching optimization scheme

Compared to other schemes, dose matching (DM) is conceptually the most closely related to optimization for greyscale response targets because it neither assumes a binary response target input nor imposes distinct algorithmic steps locally according to the response target values. The core idea of the DM formulation is to directly minimize the absolute difference between current dose response and the response target. Mathematically, the BCLP loss function is a natural generalization of DM loss function with general response mapping, local tolerance, local weighting, and global error dispersity control.

The loss function of DM takes the basic form of

$$\mathcal{L}_{DM} = \int \left| \sigma'(f(\underline{r})) - \Theta(\underline{r}) \right| d\underline{r} \tag{SE.20}$$

where $\sigma': \mathbb{R} \to \mathbb{R}$ is the material response function that maps optical dose to dose response, $\Theta: \mathbb{R}^3 \to \mathbb{R}$ is the dose response target and the integration is performed over all space. The optical tomographic dose $f(\underline{r})$ is explicitly expressed in the original work as $\frac{N_r \alpha}{\Omega} \left( T^*_{-\alpha}[g](\underline{r}) \right)$, where $T^*_{-\alpha}$ is an integral projection operator that performs the adjoint operation of the exponential Radon transform on sinogram $g$, $N_r$ is number of rotations under exposure, $\alpha$ is the absorption coefficient of the active species, and $\Omega$ is angular velocity of the rotation. The original work uses a sigmoid function $\sigma'(f) = \sigma(f - d_h, \delta) = \frac{1}{1+e^{-\frac{f-d_h}{\delta}}}$ to model nonlinear material response with parameter $d_h$ and $\delta$.

The direct correspondence of variables in $\mathcal{L}_{DM}$ in the BCLP formulation is immediately evident. $f(\underline{r})$ in both formulations represents the exact same quantity, only with $\frac{N_r \alpha}{\Omega} \left( T^*_{-\alpha}[g](\underline{r}) \right)$ simply generalized (in functional forms and in physical units) and abstracted as $P^*g$. The response target originally denoted as $\Theta$ is now denoted as $f_T$ in BCLP. The sigmoid function $\sigma'$ in DM is generalized to be a generalized logistic function $\mathcal{M}$ in BCLP. The enforcement of hard constraints on sinogram (such as non-negativity or maximum intensity) in DM can equivalently be enforced in BCLP through the definition of feasible set $S_{feasible}$ and its associated projection operation in projected gradient descent.

In the BCLP formulation, the DM loss function represents a $L_1$-norm with uniform unity weighting ($w = 1$) and zero tolerance ($\varepsilon = 0$). Under these settings, the DM loss functions can be expressed directly in the form:

$$\mathcal{L}_{DM} = \int_V \left| \mathcal{M}(f(\underline{r})) - f_T(\underline{r}) \right| d\underline{r} \, . \tag{SE.21}$$

The fact that the integral is now performed over $V = \{\underline{r} : |\mathcal{M}(f(\underline{r})) - f_T(\underline{r})| > 0\}$ instead of all space is inconsequential. The spatial locations that have zero response error ($|\mathcal{M}(f(\underline{r})) - f_T(\underline{r})| = 0$) would not contribute to the loss function or loss gradient regardless of whether they are included in the integral.

The above expressions confirmed that DM loss function coincides with a special case of BCLP. As a generalization of DM, the BCLP formulation additionally provides problem relaxation through local tolerancing and local weighting. BCLP also provides global control of error sparsity through the general form of $L_p$-norm. In the original study of DM, non-linear material response is not considered in the sinogram initialization step nor in the evaluation of chosen metrics (namely, Jaccard index and process window). The current study formally includes material response models throughout initialization (section 2.2.1), optimization, and metric evaluation (section 2.3 and S.8).

## S.11 Generalization and interpretation of the penalty minimization optimization scheme

The penalty minimization (PM) scheme optimizes dose response towards a binary target which defines the spatial regions within and outside of the printing part. The key concept in this formulation is to get the dose response to be above a certain threshold within the part and below another threshold outside of the part while ignoring the response near the part boundary. To this goal, the PM loss function applies soft constraints differently over these three distinct regions. In the original work, these three respective regions are called the eroded target object ($R_1$), eroded complement of target object ($R_2$), and buffer region. When cross-referencing with the OSMO scheme, $R_1$ and $R_2$ are similar to the in-part (IP) and out-of-part (OFP) regions in OSMO but with the regions near the part boundaries excluded.

The PM loss function only penalizes the underdosage in $R_1$ relative to a threshold $d_h$, and the overdosage in $R_2$ relative to a threshold $d_l$. In other words, $R_1$ and $R_2$ have their respective unilateral soft constraints. The PM loss function reads:

$$\mathcal{L}_{PM} = \rho_1 \int_{\sim V_1} \left(d_h - f(\underline{r})\right) d\underline{r} + \rho_2 \int_{\sim V_2} \left(f(\underline{r}) - d_l\right) d\underline{r}. \qquad (SE.22)$$

where $\sim V_1$ represents $R_1$ regions where the optical dose $f$ is lower than the threshold $d_h$ and $\sim V_2$ represents $R_2$ regions where $f$ is higher than the threshold $d_l$. $\rho_1$ and $\rho_2$ are the weightings applied on the two types of soft constraint violations. Readers should note that $\rho_1$ and $\rho_2$ are trivially renamed in this text for readability such that they correspond to violation in region $R_1$ and $R_2$ respectively.

It can be shown that the PM loss function is a special case of BCLP. This connection can be made by realizing that optimization towards a binary response target with unilateral soft constraints is a special case of optimization towards a greyscale response target with band soft constraints. As shown in the following paragraphs, all variables in PM can be equivalently reformulated to and interpreted as those in BCLP.

The volumetric optical dose $f$ in the loss function of PM carries identical meaning as DM and BCLP. As reported, the PM loss function did not apply a response mapping on $f$ such as $\sigma'(f)$ in DM or $\mathcal{M}(f)$ in BCLP. This absence of mapping is equivalent to prescribing an identity response mapping on $f$ such as $\mathcal{M}(f) = f$.

The behavior of the unilateral soft constraints in PM with thresholds $d_h$ and $d_l$ can be replicated in BCLP using the band constraint with $f_T(\underline{r})$ and $\varepsilon(\underline{r})$. In region $R_1$ where $f$ is soft constrained to be above $d_h$, the local value of $f_T$ and $\varepsilon$ in BCLP can be chosen such that the lower limit of the tolerance band would coincide with $d_h$ (i.e., $f_T - \varepsilon = d_h$). This would penalize any $R_1$ regions with $f < d_h$ based on the amount of deviation from the band limit, which can be expressed as $|f - f_T| - \varepsilon = (f_T - \varepsilon) - f = d_h - f$. The tolerance $\varepsilon$ can be chosen to be in orders of magnitude larger than the relevant scale of the response such that the upper limit of the tolerance band is never reached. $f$ (and more relevantly $\mathcal{M}(f)$) is naturally bounded (not infinite) in physical VAM settings and this allows the above construction of unilateral dose constraints from band constraints. Supplementary S.14 describes this boundedness in greater detail. Vice versa, the overdosing constraint in $R_2$ can be replicated by choosing the local value of $f_T$ and $\varepsilon$ such that the upper limit of the band coincides with $d_l$ (i.e. $f_T + \varepsilon = d_l$). In $R_2$, the value of $\varepsilon$ only needs to be large enough such that tolerance band includes minimum possible response value.

With the above spatial description of $f_T(\underline{r})$ and $\varepsilon(\underline{r})$ in place, the two disjoint constraint violation sets $\sim V_1$ and $\sim V_2$ defined in PM can be simply subsumed into one constraint violation set $V = \{\underline{r} : |\mathcal{M}(f(\underline{r})) - f_T(\underline{r})| > \varepsilon(\underline{r})\}$ in BCLP. The deviation from the tolerance band limit can represent directly the two different integrands in PM loss function:

$$|\mathcal{M}(f(\underline{r})) - f_T(\underline{r})| - \varepsilon(\underline{r}) = |f(\underline{r}) - f_T(\underline{r})| - \varepsilon(\underline{r}) = \begin{cases} f(\underline{r}) - d_l & if \ \underline{r} \in \sim V_1 \\ d_h - f(\underline{r}) & if \ \underline{r} \in \sim V_2 \end{cases}. \quad (SE.23)$$

Apart from its soft-constraint violation approach, the PM scheme further relaxes the problem by introducing buffer regions where the dose response is ignored. It chooses this buffer region to be a small neighborhood of the boundary of the binary response target such that it separates $R_1$ and $R_2$. This choice is motivated by the practical challenges in creating discontinuous dose response profiles such that the soft constraints can be met simultaneously in the whole of $R_1$ and $R_2$. This problem relaxation approach is generalized as the application of weights $w(\underline{r})$ in BCLP. In addition to the deemphasis of the buffer regions, weights $w(\underline{r})$ in BCLP also generalize the constant weights $\rho_1$ and $\rho_2$ which are applied on the soft constraint terms in $R_1$ and $R_2$ respectively. To achieve equivalent weightings in PM in BCLP, weightings should be set such that $w(\underline{r}) = \begin{cases} \rho_1, & \underline{r} \in R_1 \\ \rho_2, & \underline{r} \in R_2 \\ 0, & \underline{r} \notin (R_1 \cup R_2) \end{cases}$.

When written in BCLP variables, the whole PM loss function takes the form of a simple $L_1$-norm

$$\mathcal{L}_{PM} = \int_V w(\underline{r}) \, ||f(\underline{r}) - f_T(\underline{r})| - \varepsilon(\underline{r})| \, d\underline{r} \quad (SE.24)$$

, with $w(\underline{r})$, $f_T(\underline{r})$ and $\varepsilon(\underline{r})$ defined as in preceding paragraphs and $p = q = 1$. Hard constraints on the sinogram are enforced by the projection onto the feasible solution set $S_{feasible}$.

As reported, the PM scheme only optimizes for binary targets which have well defined boundaries. The PM scheme constructs $R_1$ and $R_2$ by eroding (a morphological operation) the binary response target and its complement, respectively, with a structuring element a few voxels in width. The buffer region is then set as the complement of the union of $R_1$ and $R_2$. Obviously, this procedural approach to designate deemphasis only applies to binary response target. One possible way to implement similar weighting relaxation with greyscale targets can be lowering the regional weights in proportion to the magnitude of the spatial gradient of the response target. In addition, regional weights can also be assigned according to the local functional importance.

## S.12 Generalization and interpretation of the Object Space Model Optimization scheme

Similar to PM, the goal of OSMO algorithm is to push the dose response of the in-part regions (denoted IP) to stay above a certain threshold $D_h$ and the dose response of the out-of-part regions (denoted OFP) to stay below a certain threshold $D_l$. The original literature did not define an explicit objective function for this algorithm. In this section, it will be shown that updates of the object space model in OSMO can also be cast into projected gradient updates of sinogram iterates for a special case of the BCLP loss function. When interpreted as a projected gradient descent, OSMO is shown to be taking a step size of 0.5. An obvious common feature that OSMO shares with PM and BCLP is that it only performs corrections on the soft constraint violation set. However, the OSMO algorithm differs from the other schemes in that it alternately handles the two constraint-violation sets (overdosing and underdosing) one at a time.

## S.12.1 Method reformulation

Here we describe the OSMO algorithm with a condensed notation, discuss the meaning of its steps, and connect it with the BCLP framework. In the following description and reformulation, singly-indexed variables ($k$ or $n$) are used in place of their doubly-indexed counterparts used in the original notation. Iteration index $k$ is used to show relationships applied to all iterations, while index $n$ is restricted to be even indices. For consistency, we count every forward propagation and backpropagation pair as one iteration, in contrast with the original article which counts two pairs as one iteration.

In the OSMO algorithm, an object space model $M_k$ of iteration $k$ continues to get updated and creates the next trial tomographic reconstruction $f_k$. The process to create $f_k$ from $M_k$ involves forward propagation, zero truncation, backpropagation and normalization. Formally, this operation is $f_k = NP^* \max(0, PM_k)$ where $N$ is a normalization operator that divides an input by its maximum, $N(a(\underline{r})) = \frac{a(\underline{r})}{\max_{\underline{r}}(a(\underline{r}))}$. Note that the model $M_k$ could contain negative values even when the binary target $f_T$ is non-negative.

The OSMO algorithm starts with $M_0 = f_T$, then proceeds to obtain $M_{n+2}$ from $M_n$ as follows:

$$f_n = NP^* \max(0, PM_n) , \tag{SE.25}$$

$$M_{n+1}(\underline{r}) = \begin{cases} (M_n(\underline{r}) - \max(0, f_n(\underline{r}) - D_l)) & \text{if } \underline{r} \in OFP \\ M_n(\underline{r}) & \text{if } \underline{r} \notin OFP \end{cases}, \tag{SE.26}$$

$$f_{n+1} = NP^* \max(0, PM_{n+1}) , \tag{SE.27}$$

$$M_{n+2}(\underline{r}) = \begin{cases} M_{n+1}(\underline{r}) + \max(0, D_h - f_{n+1}(\underline{r})) & \text{if } \underline{r} \in IP \\ M_{n+1}(\underline{r}) & \text{if } \underline{r} \notin IP \end{cases}. \tag{SE.28}$$

The algorithm continues with all $n$ being even until convergence or when satisfactory performance is reached.

The IP and OFP regions in OSMO are similar to the $R_1$ and $R_2$ regions in PM, except that the union of IP and OFP occupy all space and there is no buffer region between the two regions. From this perspective, the dose constraints $D_h$ and $D_l$ in OSMO are direct analogs of $d_h$ and $d_l$, respectively, in PM. In the BCLP formulation, they are the lower and upper bounds of the local tolerance band ($f_T \pm \varepsilon$), respectively. Although Rackson et al. [1] did not explicitly formulate OSMO as a constraint satisfaction problem as in PM, the ad-hoc penalization behavior still exists in the form of a non-linear truncation operation. By construction, the algorithm does not update regions of object space model $M$ where constraints are satisfied because the error terms $(f_n - D_l)$ or $(D_h - f_{n+1})$ would take negative values in the maximum operations $\max(0, f_n - D_l)$ and $\max(0, D_h - f_{n+1})$ respectively. This construction is equivalent to selecting only the constraint-violation regions ($\sim V_1$ and $\sim V_2$ in PM or $V$ in BCLP) in the evaluation of loss function and loss gradient.

In the following paragraphs, we will show that iterating model $M_k$ and computing its associated sinogram $g_k = \max(0, PM_k)$ is analytically equivalent to applying projected gradient descent updates directly on the sinogram iterate with an $L_2$ loss function.

To reveal the underlying gradient update operations, we change the optimization variable in the algorithm from model $M$ to sinogram $g$. We denote the update of $M_k$ as $\Delta M_k$ such that $M_{k+1} = M_k + \Delta M_k$ and then express all forward projected models $PM_k$ and $P\Delta M_k$ in terms of sinogram $g_{\pm,k} = PM_k$ and $\Delta g_{\pm,k} = P\Delta M_k$ respectively. For clarity, we denote the sinogram quantities before the enforcement of non-negativity constraint as $g_{\pm,k}$ and those after the enforcement $g_k = \max(0, PM_k)$. Now, the algorithm starts from $g_0 = \max(0, PM_0) = \max(0, Pf_T)$ and proceed with all $n$ being even:

$$f_n = NP^* g_0 \qquad (SE.29)$$

$$\Delta M_n(\underline{r}) = \begin{cases} -\max(0, f_n(\underline{r}) - D_l) & \text{if } \underline{r} \in OFP \\ 0 & \text{if } \underline{r} \notin OFP \end{cases} \qquad (SE.30)$$

$$g_{n+1} = \max(0, PM_n + P\Delta M_n) = \max(0, g_{\pm,n} + \Delta g_{\pm,n}) \qquad (SE.31)$$

$$f_{n+1} = NP^* g_{n+1} \qquad (SE.32)$$

$$\Delta M_{n+1}(\underline{r}) = \begin{cases} \max(0, D_h - f_{n+1}(\underline{r})) & \text{if } \underline{r} \in IP \\ 0 & \text{if } \underline{r} \notin IP \end{cases} \qquad (SE.33)$$

$$g_{n+2} = \max(0, PM_{n+1} + P\Delta M_{n+1}) = \max(0, g_{\pm,n+1} + \Delta g_{\pm,n+1}) \qquad (SE.34)$$

In this rewritten form, it is easy to see that sinogram $g$ is progressively updated to correct for $IP$ and $OFP$ violation sets alternately. The magnitude of the correction is proportional to the dose response error from the thresholds and has the form of $D_l - f_k$ (negative sign incorporated) in OFP or $D_h - f_{k+1}$ in IP.

## S.12.2 Equivalent quantities in the BCLP formulation

The above operations can be similarly achieved by projected gradient descent updates with a special case of BCLP loss function $\mathcal{L}$. This special case of $\mathcal{L}$ would have $p = q = 2$, and $\mathcal{M}(f) = f$. Identical to the discussion on PM, $f_T$ and $\varepsilon$ are declared such that $f_T(\underline{r}) + \varepsilon(\underline{r}) = D_l$ for $\underline{r} \in OFP$ and $f_T(\underline{r}) - \varepsilon(\underline{r}) = D_h$ for $\underline{r} \in IP$. The tolerance $\varepsilon$ is chosen to be large enough that it is numerically improbable to have OFP underdose and IP overdose violations (as described in supplementary S.14). This setting of $f_T$ and $\varepsilon$ effectively reproduces the unilateral soft constraints in OSMO.

With the above parameters, the BCLP loss function and its gradient (derived in supplementary S.1) at iteration $k$ are written as

$$\mathcal{L}_k = \int_{V_k} w(\underline{r}) \big| |f_k(\underline{r}) - f_T(\underline{r})| - \varepsilon(\underline{r}) \big|^2 d\underline{r} \qquad (SE.35)$$

$$\left(\nabla_g \mathcal{L}(\underline{r}')\right)_k = 2P\left(v_k(\underline{r}) w(\underline{r}) (|f_k(\underline{r}) - f_T(\underline{r})| - \varepsilon(\underline{r})) \, sgn(f_k(\underline{r}) - f_T(\underline{r}))\right)(\underline{r}') \qquad (SE.36)$$

where $v_k(\underline{r}) = \begin{cases} 1 & \text{if } \underline{r} \in V_k \\ 0 & \text{if } \underline{r} \notin V_k \end{cases}$.

Applying projected gradient method with step size $\eta$ and $S_{feasible} = \{g \in S \mid g(\underline{r}') \geq 0 \; \forall \underline{r}'\}$ would give the following update for iteration $k$:

$$g_{k+1} = S_{feasible}\left(g_k - \eta \left(\nabla_g \mathcal{L}\right)_k\right)$$

$$= \max\left(0, g_k - \eta \, 2P\left(v_k(\underline{r})w(\underline{r})\left(|f_k(\underline{r}) - f_T(\underline{r})| - \varepsilon(\underline{r})\right) sgn\left(f_k(\underline{r}) - f_T(\underline{r})\right)\right)\right)$$

$$= \max\left(0, g_k + P\left(2\eta v_k(\underline{r})w(\underline{r})\left(f_T(\underline{r}) - f_k(\underline{r}) + \varepsilon(\underline{r})sgn\left(f_k(\underline{r}) - f_T(\underline{r})\right)\right)\right)\right) \quad (SE.37)$$

As mentioned earlier in the section, only regions with soft constraint violated would influence the update of optimization variables. Recalling $f_T(\underline{r}) + \varepsilon(\underline{r}) = D_l$ for $\underline{r} \in OFP$ and $f_T(\underline{r}) - \varepsilon(\underline{r}) = D_h$ for $\underline{r} \in IP$ allows us to see how the soft band constraint of BCLP provides the correction terms in OSMO:

$$f_T(\underline{r}) - f_k(\underline{r}) + \varepsilon(\underline{r})sgn\left(f_k(\underline{r}) - f_T(\underline{r})\right) = \begin{cases} D_l - f_k(\underline{r}) < 0 & \text{for } \underline{r} \in (V_k \cap OFP) \\ D_h - f_k(\underline{r}) > 0 & \text{for } \underline{r} \in (V_k \cap IP) \end{cases}. \quad (SE.38)$$

With the choice of $\eta = 1/2$ and $w(\underline{r}) = 1$,

$$g_{k+1} = \max\left(0, g_k + P\left(v_{k,OFP}(\underline{r})\left(D_l - f_k(\underline{r})\right)\right) + P\left(v_{k,IP}(\underline{r})\left(D_h - f_k(\underline{r})\right)\right)\right) \quad (SE.39)$$

, where $v_{k,OFP}(\underline{r}) = \begin{cases} 1 & \text{if } \underline{r} \in (V_k \cap OFP) \\ 0 & \text{otherwise} \end{cases}$ and $v_{k,IP}(\underline{r}) = \begin{cases} 1 & \text{if } \underline{r} \in (V_k \cap IP) \\ 0 & \text{otherwise} \end{cases}$.

In eq. (SE.39), $P\left(v_{k,OFP}(\underline{r})\left(D_l - f_k(\underline{r})\right)\right)$ and $P\left(v_{k,IP}(\underline{r})\left(D_h - f_k(\underline{r})\right)\right)$ are update terms of $g_{k+1}$ that correspond to the correction for constraint violation in OFP and IP, respectively. These two update terms are applied in every iteration.

Comparing the above special case of BCLP and the rewritten form of OSMO, we find that the correction terms in OSMO can be generated by projected gradient descent on an implicit $L_2$ loss function. In particular, the backpropagation $P^*$ and forward propagation $P$ steps in OSMO occur in the evaluation of the dose tomogram $f$ and the gradient $\nabla_g \mathcal{L}$, respectively. The model update $M_{k+1} = M_k + 1(\Delta M_k)$ implied that the descent step size $\eta$ is ½. In BCLP, the evaluation of $f$ does not require a normalization step and preserves the physical dose unit (supplementary S.2).

As stated, the projected gradient update on eq. (SE.39) differ from OSMO in that it handles the two types of soft constraint violation (OFP overdose and IP underdose) together in every iteration instead of separately in even and odd iteration in OSMO (eq.(SE.31) and eq.(SE.34)). Users can opt to reproduce this behavior in BCLP by extending the definition of weight $w(\underline{r})$ to be a function $w(\underline{r}, k) = w_k(\underline{r})$ of iteration number $k$.

The definition of $w_k(\underline{r}) = \begin{cases} 1 & \text{if } (k \text{ is even}) \text{ and } (\underline{r} \in OFP) \\ 1 & \text{if } (k \text{ is odd}) \text{ and } (\underline{r} \in IP) \\ 0 & \text{otherwise} \end{cases}$ allows this alternate handling such that:

$$g_{n+1} = \max\left(0, g_n + P\left(v_{n,OFP}(\underline{r})\left(D_l - f_n(\underline{r})\right)\right)\right) = \max(0, g_n + \Delta g_{\pm,n}) \,, \quad (SE.40)$$

$$g_{n+2} = \max\left(0, g_{n+1} + P\left(v_{n+1,IP}(\underline{r})\left(D_h - f_{n+1}(\underline{r})\right)\right)\right) = \max(0, g_{n+1} + \Delta g_{\pm,n+1}) \ . \quad (SE.41)$$

With alternate weighting $w_k$, BCLP (eq.(SE.40) and (SE.41)) employs an almost identical updating step as in OSMO (eq.(SE.31) and (SE.34)). Both schemes produce the same gradient update terms ($\Delta g_{\pm,n}$ and $\Delta g_{\pm,n+1}$). However, BCLP (and other schemes such as DM and PM) always add the gradient update $\Delta g_{\pm,n}$ to an iterate $g_n$ that is inside the feasible region $S_{feasible}$ while OSMO adds the gradient update $\Delta g_{\pm,n}$ to an iterate $g_{\pm,n}$ that may or may not be inside the feasible region $S_{feasible}$. This subtle difference arises from the fact that the OSMO algorithm is designed to store the optimization iterate $M_k$ in tomogram function space ($T$) instead of sinogram function space ($S$) where the feasible set is defined. Overall, OSMO still minimizes the above special case of the BCLP loss function by traversing solutions in the direction of steepest descent until the loss gradient vanishes.

Supplementary S.13 compares BCLP optimization runs with and without alternate handling. The comparison shows no obvious benefit of implementing alternate handling in the test case.

## S.13 Effect of alternate handling of positive and negative error

The OSMO optimization scheme handles the underdosing error in IP and overdosing error in OFP alternately. Although the BCLP framework can reproduce this behavior by setting an iteration-dependent weight $w(\underline{r})$ (as discussed at the end of section S.12.2), it is unclear whether this alternate handling is beneficial to convergence. This section compares the optimization convergence behavior with and without implementing this alternate handling.

For relevance to application context of OSMO, this comparative study run two BCLP optimizations with settings that matches those in OSMO ($p = q = 2$, step size $\eta = 1/2$, linear material response, and unilateral dose constraints). In the first run, weight is set to be one everywhere such that both errors in IP and OFP are handled in every iteration. In the second run, weight is set to be

$$w_k(\underline{r}) = \begin{cases} 1 & if\ (k\ is\ even)\ and\ (\underline{r} \in OFP) \\ 1 & if\ (k\ is\ odd)\ and\ (\underline{r} \in IP) \\ 0 & otherwise \end{cases}$$

for iteration $k$. This iteration-dependent weight implies that OFP overdosing error and IP underdosing error are handled in even and odd iterations, respectively. The response target comprises four binary gratings with local values being {0,1}. The material response model $\mathcal{M}$ is an identity function of dose $f$. Both runs terminate at the 250-th iteration with no other imposed termination criterion.

These two runs use a common initial solution $g_0$ that is generated through steps detailed in section 2.2.1. The initialization computation used the raw response target mentioned above with values in {0,1}.

For both optimization runs, the target values $f_T$ and local tolerance are set such that the band constraints in BCLP framework effectively reproduce the unilateral constraint in OSMO scheme. By setting $f_T = 10^6$ and $\varepsilon = 10^6 - 0.8$ in the IP region, the lower limit of the band is $f_T - \varepsilon = 0.8$ and the upper limit of the band is $f_T + \varepsilon = 2 * 10^6 - 0.8 \approx 2 * 10^6$ in IP. Effectively, only the lower tolerance limit is active in IP during optimization. For OFP, since the dose and dose response never go negative, $f_T = 0$ and $\varepsilon = 0.2$ provides a unilateral constraint with upper limit at 0.2. In short, the IP region only has a lower tolerance limit of 0.8 and the OFP region only has an upper tolerance limit of 0.2.

Supplementary Fig. 1 shows the binary target, the response of the two runs and the convergence plot of the two runs. On the figure, the first run with constant weight is labelled as CW run while the second run with alternating weight is labelled as AW run.

Since the weight alternates from iteration to iteration, the loss value in the AW run also jumps when measuring error in IP and OFP alternately. Initially, the iterates in the AW run have no IP underdose error so the odd iterations have a zero loss. As optimization progresses, both odd and even iterations have non-zero losses as IP underdosing starts to occur. In the AW run, an $L_2$ evaluation metric is used to replicate the CW loss function and evaluates the AW run solutions with constant weight. This metric provides values that can be directly compared to the loss in the CW run. Although the AW run has a fluctuating loss value during optimization, its evaluation metric showed a relatively smooth decrease.

Comparing the $L_2$ evaluation metric of AW run (in blue) and loss value of CW (in green), the CW run converged almost twice as fast. This can be explained by the fact that the loss gradient in the CW run corrects for both types of error in every iteration while AW run corrects one type of error in each iteration. Quantitatively, the loss gradient in CW run (as expressed in eq. (SE. 39)) includes more corrective information per iteration than that of the AW run (as in eq.(SE. 40) and (SE. 41)). In this test, we observe no obvious benefit in handling overdosing and underdosing error alternately.

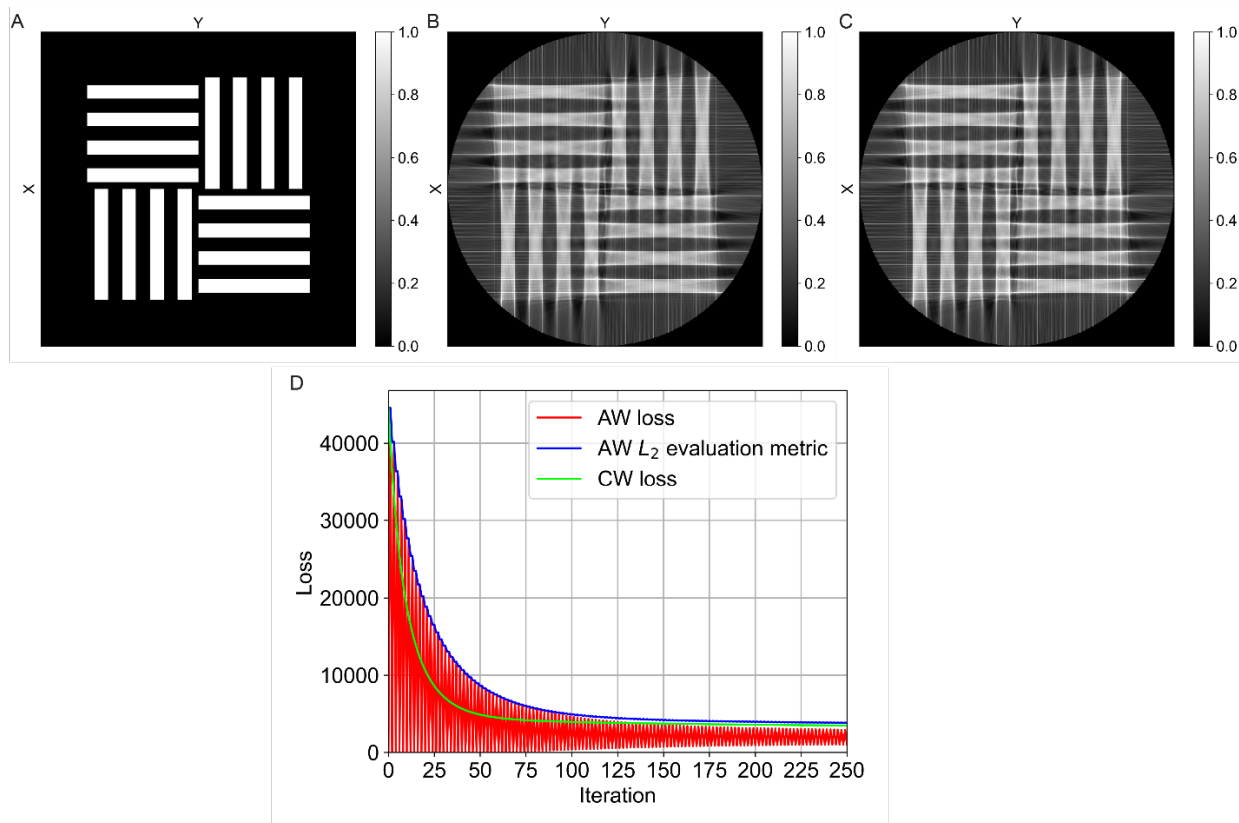

**Supplementary Fig. 1.** A) Response target, B) Response of the solution at 250-th iteration for the CW run, C) Response of the solution at 250-th iteration for the AW run, D) Convergence plot of the two runs. The loss value for CW run and AW run is colored green and red, respectively. For fair comparison, a $L_2$ evaluation metric (colored blue) replicates the loss function of the CW run and evaluates the solution in the AW run.

## S.14 Boundedness of $\mathcal{M}(f)$ for implementation of soft unilateral dose constraints

In both the discussion of PM and OSMO (in section S.11 and S.12), we argue that the behavior of unilateral dose constraints can be achieved by certain choices of tolerance band with variables $f_T$ and $\varepsilon$. To completely replicate the unilateral dose constraints in $R_1$, the upper limit of the tolerance band must be assigned such that it is larger than any possible value of $\mathcal{M}(f)$. This in turn requires $\mathcal{M}(f)$ be bounded (not infinite). Here, we delineate the reason why this requirement is always satisfied in physical systems.

$\mathcal{M}(f)$ models the physical response of the material to exposure dose $f$. As an output of the function $\mathcal{M}(f)$, the physical response such as degree-of-conversion (DOC), elastic modulus, and refractive index always has an upper bound. For example, DOC has a maximum value of 100% by definition. The existence of such upper bounds motivates the use of a response model with saturation, such as a logistic function. This work uses a generalized logistic function (S.6) as the response model, where the right asymptote (parameter $K$) designates the upper bound of the response value. When the local tolerance band covers this upper bound, the optimization behavior of the band constraint at that location is identical to a unilateral constraint. Even for a hypothetical material where it has an unbounded response that varies monotonically with $f$, the boundedness of input dose $f$ still bounds the output response value $\mathcal{M}(f)$. The dose $f$ is always bounded due to finite available optical power and fabrication time.

Conversely, it is obvious that the unilateral dose constraints in $R_2$ can also be reproduced by setting a tolerance band that covers the lower bound of response.

## S.15 Convergence plot of parameter sweep on steepness of material response

Supplementary Fig. 2 plots the history of loss value for all optimization runs in the sweep of parameter $B$ in the material response. The run with $B = 150$ significantly lowers the loss within the first 1000 iterations but fails to meet the convergence criterion in 2000 iterations.

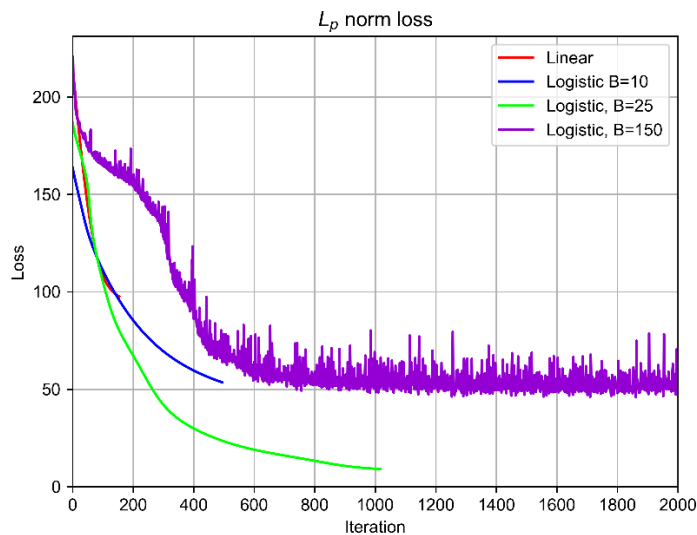

**Supplementary Fig. 2.** Convergence plot of all runs in sweep of steepness of material response.

## S.16 Optimization Demonstration for a 3D binary target

A demonstrative optimization run was performed for a 3D binary response target at 240×240×256 resolution. The optimization used the default material response parameters as listed in Supplementary Table 1, except that $B = 25$. Other optimization parameters were: $p = 2$, $q = 1$, $\varepsilon = 0.1$, $w = 1$, and $\eta(step\ size) = 3000$. The loss trajectory did not meet the convergence criterion before the optimization terminated at the 50-th iteration with a final loss value of 227.5. Supplementary Fig. 3 shows the binary target and the resulting response from the optimization. The response plot on the figure is visualized with linear opacity from response 0 to 1.

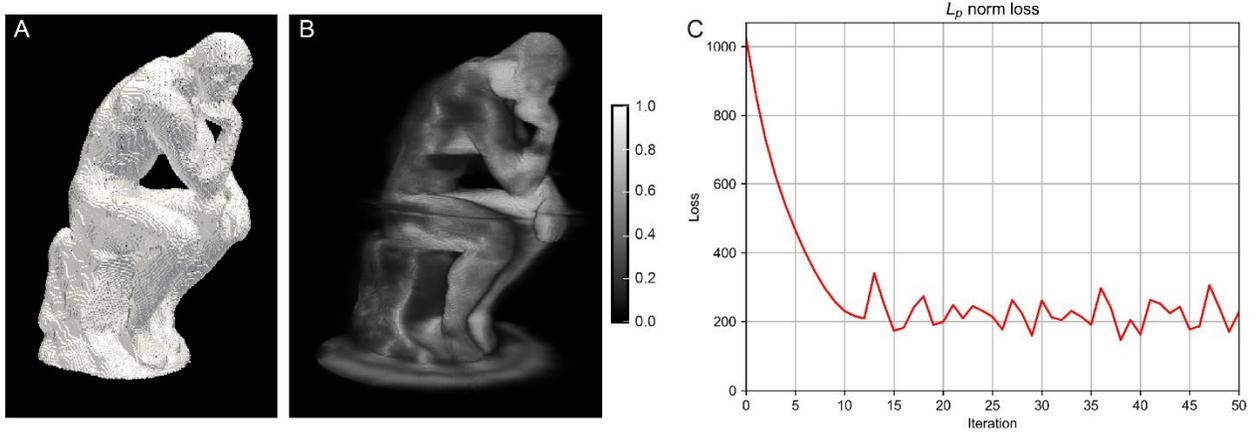

**Supplementary Fig. 3.** Demonstration of optimization for a binary 3D target. (A) Binary response target, (B) Reconstructed response, (C) Plot of loss value at each iteration. The reconstructed response in (B) is colorized according to the response value and displayed with a linear opacity increasing from 0 to 1 in the range of response. Both (A) and (B) are visualized in the visualization software tomviz.

## S.17 Alternative interpretation of the BCLP loss function

As expressed in eq. (1), the loss function can be interpreted as a weighted $L_p$-norm of an error term, where only the error beyond a tolerance is counted. The loss function expression is a direct generalization of the loss function in DM and PM. As a volume integral, this expression emphasizes the spatial dependence of each quantity.

Here we provide an alternative expression and interpretation that emphazise the statistical nature of the loss function. For simplicity, we seek an alternate expression of $\mathcal{L}^{\frac{p}{q}}$ instead of $\mathcal{L}$. The monotonicity of $\mathcal{L}^{\frac{p}{q}}$ with respect to $\mathcal{L}$ (for all $\frac{p}{q} > 0$) implies that $\mathcal{L}^{\frac{p}{q}}$ and $\mathcal{L}$ share the same locations of local and global minima.

First, let $E$ be the error term $E(\mathbf{r}) = |\mathcal{M}(f(\mathbf{r})) - f_T(\mathbf{r})| - \varepsilon(\mathbf{r})$. $E$ is always non-negative when $\varepsilon(\mathbf{r}) = 0$ but can be negative with $\varepsilon(\mathbf{r}) > 0$.

Then, we define a function $h$ to measure the weighted volume in all space that takes a specific error value $E'$. In other words, the function measures a total weighted volume of the level set of $E$.

$$h(E') = \int_\infty \delta(E(\mathbf{r}) - E')\, w(\mathbf{r}) d\mathbf{r} \geq 0, \qquad (SE.42)$$

where $\delta$ is the Dirac delta function. In discrete form, this function is equivalent to a histogram of $E$ where the differential volume $d\underline{r}$ is locally scaled by weight $w(\underline{r})$. Readers should note that only the part of the histogram on positive $E'$ is relevant because the population with $E' \leq 0$ is within tolerance.

Next, we introduce a relevant concept of moments which is widely used in statistics. Moment in mathematics is similar to moment in mechanics, where a function is multiplied with the coordinate. The $p$-th moment of a function $a(x)$ is $\int_\infty a(x)x^p \, dx$. If $a(x)$ is mass density in mechanics, then the zeroth, first and second moment correspond to the total mass, center of mass (times total mass), and mass moment of inertia, respectively. If $a(x)$ is probability density, the zeroth, first and second moments correspond to the total probabilty (which always equals one), mean, and variance, respectively. Typically, the variance is measured about the mean instead of the origin and is the second central moment.

With this concept, we can express the positive $p$-th moment of the histogram of $E$ as:

$$M_p = \int_0^\infty h(E')E'^p \, dE'. \qquad (SE.43)$$

It can be shown that this positive $p$-th moment is equivalent to $\mathcal{L}^{\frac{p}{q}}$ which the optimization minimizes.

$$M_p = \int_0^\infty E'^p h(E') \, dE' = \int_0^\infty E'^p \int_\infty \delta(E(\underline{r}) - E') w(\underline{r}) d\underline{r} \, dE'. \qquad (SE.44)$$

$E'$ is not a function of $\underline{r}$ and can be moved inside the inner integral; therefore

$$M_p = \int_0^\infty \int_\infty E'^p \delta(E(\underline{r}) - E') w(\underline{r}) \, d\underline{r} \, dE'. \qquad (SE.45)$$

If we change the order of integration,

$$M_p = \int_\infty \left( \int_0^\infty E'^p \delta(E(\underline{r}) - E') \, dE' \right) w(\underline{r}) d\underline{r}. \qquad (SE.46)$$

Then by the sifting property of the delta function,

$$M_p = \int_\infty v(E(\underline{r})) E(\underline{r})^p w(\underline{r}) \, d\underline{r} = \int_V w(\underline{r}) |E(\underline{r})|^p \, d\underline{r}, \qquad (SE.47)$$

where $v$ is the indicator function $v(E(\underline{r})) = \begin{cases} 1 \text{ if } E(\underline{r}) > 0 \\ 0 \text{ if } E(\underline{r}) \leq 0 \end{cases}$ and $V = \{\underline{r} : E(\underline{r}) > 0\}$. The indicator function appears because the above sifting is performed on an integration with integration limits from $E' = 0$ to $E' \to \infty$ instead of extending to infinities on both sides.

From the last expression of $M_p$, we can recognize that:

$$M_p = \int_V w(\underline{r}) |E(\underline{r})|^p \, d\underline{r} = \mathcal{L}^{\frac{p}{q}}. \qquad (SE.48)$$

Therefore, minimization of the loss function is minimizing the positive $p$-th moment of the weighted histogram of the error term $E$. As the optimization progresses, the positive error population moves to zero or negative $E$ on the histogram.

This alternative expression provides an interesting interpretation of the loss function parameters. The weight $w(\mathbf{r})$ scales the local differential volume and hence the local contribution to the histogram $h$. The positive and negative parts of the histogram represent the population outside and inside of the tolerance band, respectively. With a zero tolerance $\varepsilon$, the entire population resides on $E \geq 0$ on the histogram. A greater tolerance $\varepsilon$ makes $E$ more negative. On the histogram, a greater $\varepsilon$ shifts the entire error population of $E$ to the negative side, which does not contribute to the positive moment. Finally, the greater the $p$ value, the greater the moment contributed by the large-error population. Therefore, a greater $p$ value puts more emphasis on the population with large errors. These interpretations are consistent with the findings in the parameter sweep. The alternative expression in terms of moment abstracts away the spatial dependence and provides a statistical perspective to visualize the physical meaning of the BCLP loss function.

## S.18 Sinogram discretization and hard constraints modeling

This section provides some general guidelines on selecting proper sinogram discretization. In particular, these guidelines are made for the most commonly used parallel-beam configuration. Although the demonstration of this work only applied the simplest hard constraint (the non-negativity constraint) on sinogram variables, practical situations also impose other hard constraints. The second half of this section also discusses how the commonly encountered hard constraints can implemented in this optimization framework.

### S.18.1 Transverse and angular discretization of sinogram in parallel-beam configuration

The transverse and angular discretization of the sinogram heavily depends on the spatial discretization of the tomogram domain which houses the response target. In other to completely reconstruct the response target at the specified resolution, the Nyquist sampling criterion must be satisfied. The following paragraphs describe this resolution requirement for the case of parallel-beam configuration under the assumption that both the voxel size and the number of voxels of the tomogram domain are isotropic.

Let the tomogram domain be discretized in $N \times N \times N$ voxel grid that spans $L \times L \times L$ length. The original discretization of the tomogram is correspondingly $L/N$. If a parallel projection of the sinogram has to sample this signal without loss of information in its transverse direction, the transverse discretization of the sinogram must be at least as fine as that of the tomogram ($L/N$). In other words, if the width of the projection is also $L$ (which is the common practice), the number of pixels on the sinogram should be at least $N$.

On the other hand, the required angular discretization of the sinogram is established by the sampling requirement in the Fourier domain. Given by the Fourier slice theorem, each projection is in fact encoding information that falls on a slice in the Fourier space of the tomogram. Each projection oriented at a particular angle would correspond to a slice in Fourier space that is oriented at the same angle and centered at the origin. Therefore, the angular discretization of the sinogram is dictating the angular discretization at which the Fourier space of the tomogram is sampled. If the angular discretization is too low, certain spatial frequency information in the Fourier space is discarded, and this phenomenon is called aliasing.

More concretely, the number of required angular samples is also found to be proportional to $N$. By the properties of discrete Fourier transform, the number of voxels of the Fourier space of the tomogram is also $N \times N \times N$. If the transverse discretization of the sinogram is high enough, the slices that the projections correspond to would reach the end of this Fourier domain. For the sampled slice to be dense enough to not miss any of these voxels, the slices must have an angular separation less than $\approx 2/N$ radians. Since the slices are centered at the origin and reach both sides of the domain, the slices only need to sample over 180° ($\pi$ radians) to completely cover the Fourier domain. In other words, the total number of required angular positions would be $\pi N/2$. In summary, both the number of transverse and angular samples are directly proportional to the number of voxels across the tomogram.

This sampling requirement is more rigorously discussed in the context of CT[35,36]. The original VAM publication also discussed this sampling requirement[26].

### S.18.2 Hard constraints on the set of feasible solutions

In this work, the set of feasible solution $S_{feasible}$ is only confined to the set of non-negative sinogram $\{g \in S \mid g(\underline{r}') \geq 0 \ \forall \underline{r}'\}$. In practice, the capability of the available projection hardware would also limit the range and bit-depth of the achievable sinogram value. Mathematically, these hard constraints would further confine the set $S_{feasible}$.

If the projection has a minimum dark level and hence delivers a minimum dose of $h_{min}$, the corresponding feasible set would be $\{g \in S \mid g(\underline{r}') \geq h_{min} \ \forall \underline{r}'\}$. In addition, if the maximum dose is also limited by the finite power of the projector and finite exposure time to be $h_{max}$, the feasible set would read as $\{g \in S \mid h_{max} \geq g(\underline{r}') \geq h_{min} \ \forall \underline{r}'\}$. In the step where the optimization algorithm projects the latest sinogram onto $S_{feasible}$, the computer program should clip the out-of-bounds values accordingly.

If the projector has a finite bit-depth $b$, the dose level should also be discretized into $2^b$ levels. The corresponding feasible set of sinograms would be $\{g \in S \mid h_{max} \geq g(\underline{r}') \geq h_{min}, \ g(\underline{r}') \in \left\{\frac{c}{2^b}(h_{max} - h_{min}) + h_{min} \mid c \in \mathbb{Z}_{\geq 0}\right\} \forall \underline{r}'\}$. In programming, projecting onto these levels is performed by a rounding operation that use $(h_{max} - h_{min})/2^b$ as the unit. It should be noted that this particular constraint turns the optimization into an integer programming problem and hinders gradient descent. If the optimization is found to get stuck and converge prematurely, it is commonly recommended to only apply this constraint at the end of the optimization.

### S.19 Methods to account for change in refractive index in optical propagation

When material polymerizes during light exposure, the change in refractive index may refract light appreciably and reduce the accuracy of dose delivery at the later stage of the patterning process. This is particularly important for greyscale printing and for materials with high refractive index changes.

Although the demonstrations presented in this paper assumed a static propagation model, the proposed optimization framework can also be used to generate projections for dynamic optical media, as long as the backpropagation operator $P^*$ is updated according to the refractive index changes. In order to perform such an update, the refractive index has to be determined either by measurement or simulation. The following paragraphs discuss some of the potential dynamic compensation methods based on measured or simulated refractive index changes.

### S.19.1 Approaches based on refractive index measurement

With appropriate hardware, it is possible to measure the local changes in refractive index to provide real-time process feedback. Since this measurement has to be performed for locations that are buried deep in the volume, traditional contact-based refractometry methods are not applicable. Fortunately, there are tools that are developed for volumetric measurements of refractive index changes. Prime examples of these tools are Color Schlieren Tomography (CST)[37,38], Optical Scattering Tomography (OST)[39,40], and Optical Coherence Refraction Tomography (OCRT)[41]. CST and OST detect the refractive index changes in the patterning volume through Schlieren imaging and dark-field scatterometry, respectively. By taking measurements from different angles relative to the material container and using these data to perform tomographic reconstruction, these methods can directly reconstruct the change of refractive index at every point within the volume and continuously in time. In comparison, OCRT is an iterative method. It first takes optical coherence tomography (OCT) measurements from different angles relative to the material container, and then deduces the refractive index distribution by iteratively registering the simulated OCT image with the measured OCT image. Since all these measurements are naturally performed from multiple angles similar to the projection of tomographic VAM, they can be run continuously and in parallel to the printing process as the material rotates.

With the incoming feed of refractive index measurement, an optical simulation can then be performed to update the backpropagation operator $P^*$. Typically, the photopolymerization in VAM would create refractive index gradients that bend the trajectory of the passing rays into curves. Therefore, the optical simulation must be able to model these refraction effects caused by index gradients. One of the prominent methods to perform this simulation is Eikonal ray tracing[42,43], which computes the ray trajectory by solving the underlying ordinary differential equations. In the context of VAM, the authors of this work have previously demonstrated the potential of Eikonal ray tracing to account for attenuation and refraction phenomena in projection optimization[44]. Following similar procedures, users can program the computer to perform this optimization in real time as the print progresses and as new index measurements are made.

Certainly, this measurement-based approach would require the refractive index measurement, optical simulation, and projection optimization to occur at a rate that is relevant to the polymerization process itself. Currently, this heavy computational requirement remains a technical barrier to achieving closed-loop VAM. Nevertheless, these challenges can likely be overcome by refinement of measurement techniques, algorithmic improvements and the use of parallel processing.

## S.19.2 Approaches based on reaction simulation and modeling

Another way to estimate the refractive index is by material characterization and simulation of the polymerization reaction. Similar to the measurement-based approach, the projection optimization and reaction simulation could take place alternately and step through time progressively to account for the latest changes. One major advantage of simulating the change instead of measuring it is that the entire process planning step can be done asynchronously ahead of time and hence does not require real-time computing.

The foundation of the simulation can be based on either an empirical or a mechanistic model. The empirical model can be built by directly characterizing the evolution of the refractive index due to light exposure through photorefractometry. In comparison, the mechanistic model estimates such changes by modeling the photopolymerization kinetics and correlating refractive index with the change in driving quantities such as degree-of-conversion (DOC) and temperature. Evidently, the mechanistic model is substantially more complex to build and run. Nevertheless, methods for such high-fidelity reaction simulations have already been developed for the VAM[22].

Overall, the full computational procedure can be listed as follows. At the beginning, an initial run of the optimization provides an optimized projection set to the process simulation. In turn, the process simulation is run for a selected time step and provides an updated prediction of refractive index. This index distribution can be used to update the backpropagation operator through optical simulation (such as Eikonal ray tracing mentioned above) to account for the refraction effects. Then a following projection optimization run is conducted using this new operator while accounting for the history of previously delivered dose. If the DOC is also provided by the simulation, the DOC distribution can directly be converted to a starting response value in the optimization to account for previous exposure (in place of dose history). After certain iterations of optimization, the refined projection set can be used to advance the simulation by another time step. By performing simulation and optimization alternately until, a sequence of optimized projection sets is obtained, which can be displayed in the actual printing process. Supplementary Fig. 4 summarizes the overall computational procedure with measurement-based and simulation-based compensation approaches.

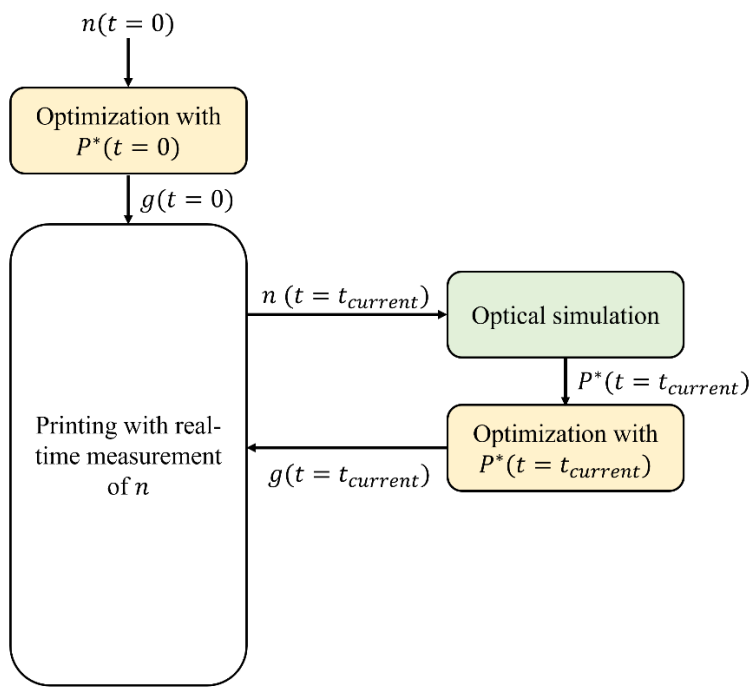
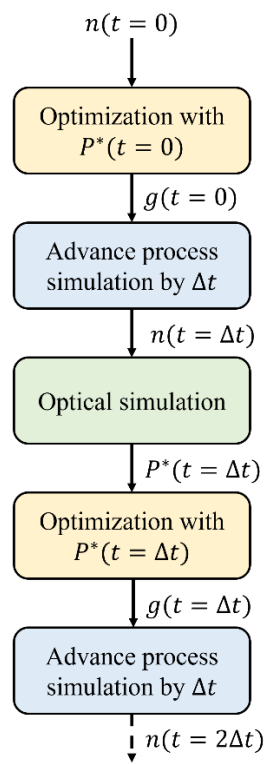

**Supplementary Fig. 4.** A flow chart of computational steps to account for refractive index change via measurement and simulation approach. In the figure, the refractive index distribution is denoted as $n$, time is denoted as $t$, backpropagation operator is denoted as $P^*$, and the optimized projection set at a given time is denoted as $g$.

Obviously, the simulation approach also has its technical limitations. Currently, high-fidelity simulation such as VirtualVAM[22] is performed on a high-performance computing (HPC) cluster and requires tens of hours of computational time. This high cost limits the practicality of mechanistic models. Therefore, the empirical or reduced-order models remain the practical options for the general public in terms of computational cost.

Lastly, it should be noted that it is possible to have a mix of both measurement and simulation in the real-time feedback loop. This hybrid approach can leverage the robustness of the measurement approach and the predictive power of simulation to estimate the upcoming trajectory of the refractive index and make corrections accordingly.